\documentclass[reqno, 11pt]{amsart}
\usepackage{amssymb}

\usepackage[pdftex]{graphicx}
\usepackage{listings}
\usepackage{multirow}
\usepackage{placeins}
\usepackage{color}
\usepackage{lscape}
\usepackage{dsfont}

\usepackage{graphicx}
\usepackage{subcaption}

\usepackage{lipsum}

\newcommand{\BlackBoxes}{\global\overfullrule5pt}

\BlackBoxes

\newcommand{\R}{\mathbb{R}}

\newcommand{\abs}[1]{\left| #1 \right|}

\usepackage{bm}

\newcommand{\vb}[1]{\bm{#1}}

\newcommand{\diag}{\operatorname{diag}}  
\newcommand{\expo}{\operatorname{exp}}

\newcommand{\Q}{{\mathbb Q}}

\newtheorem{thm}{Theorem}[section]

\newtheorem{lem}[thm]{Lemma}

\theoremstyle{definition}
\newtheorem{de}[thm]{Definition}

\newtheorem{rem}[thm]{Remark}

\numberwithin{equation}{section}

\def\0{\kern0pt\-\nobreak\hskip0pt\relax}

\makeatletter
\AtBeginDocument{%
 \def\@serieslogo{%
 \vbox to\headheight{%
 \parindent\z@ \fontsize{6}{7\p@}\selectfont
 \vss}}}

\def\makeoverbar#1#2#3#4#5#6#7{%
 \setbox0=\hbox{$\m@th#2\mkern#5mu{{}#3{}}\mkern#6mu$}%
 \setbox1=\null \dimen@=#4\fontdimen8#13 \dimen@=3.5\dimen@
 \advance\dimen@ by \ht0 \dimen@=-#7\dimen@ \advance\dimen@ by \wd0
 \ht1=\ht0 \dp1=\dp0 \wd1=\dimen@
 \dimen@=\fontdimen8#13 \fontdimen8#13=#4\fontdimen8#13
 \rlap{\hbox to \wd0{$\m@th\hss#2{\overline{\box1}}\mkern#5mu$}}
 \fontdimen8#13=\dimen@}

\def\mylabel#1#2{{\def\@currentlabel{#2}\label{#1}}}

\makeatother

\newcommand{\vecop}{\operatorname{vec}}

\DeclareMathOperator{\EX}{\mathbb{E}}

\allowdisplaybreaks

\begin{document}


\title[Utility maximization in multivariate Volterra Models]{Utility maximization in multivariate Volterra Models}

\author[F. \smash{Aichinger}]{Florian Aichinger${}^*$}
\address[F. Aichinger]{Institute for Financial Mathematics and Applied Number Theory, University of Linz, AT-4040 Linz, Austria}

\email{florian.aichinger@jku.at}

\author[S. \smash{Desmettre}]{Sascha Desmettre${}^\dagger$}
\address[S. Desmettre]{Institute for Financial Mathematics and Applied Number Theory, University of Linz, AT-4040 Linz, Austria}

\email{sascha.desmettre@jku.at}

\thanks{${}^*$ Institute for Financial Mathematics and Applied Number Theory, University of Linz, AT-4040 Linz, Austria}
\thanks{${}^\dagger$ Institute for Financial Mathematics and Applied Number Theory, University of Linz, AT-4040 Linz, Austria}

\thanks{}

\date{\today}

\begin{abstract}
This paper is concerned with portfolio selection for an investor with power utility in multi-asset financial markets in a rough stochastic environment. We investigate Merton’s portfolio problem for different multivariate Volterra models, covering the rough Heston model. First we consider a class of multivariate affine Volterra models introduced in [E. Abi Jaber et al., SIAM J. Financial Math., 12, 369–409, (2021)]. Based on the classical Wishart model described in [N. B\"auerle and Li, Z., J. Appl. Probab., 50, 1025–1043 (2013)], we then introduce a new matrix-valued stochastic volatility model, where the volatility is driven by a Volterra-Wishart process. Due to the non-Markovianity of the underlying processes, the classical stochastic control approach cannot be applied in these settings. To overcome this issue, we provide a verification argument using calculus of convolutions and resolvents. The resulting optimal strategy can then be expressed explicitly in terms of the solution of a multivariate Riccati-Volterra equation. We thus extend the results obtained by Han and Wong to the multivariate case, avoiding restrictions on the correlation structure linked to the martingale distortion transformation used in [B. Han and Wong, H. Y., Finance Res. Lett., 39 (2021)]. We also provide existence and uniqueness theorems for the occurring Volterra processes and illustrate our results with a numerical study.
\end{abstract}

\maketitle

\vspace{0.5cm}
\begin{minipage}{14cm}
{\small
\begin{description}
\item[\rm \textsc{ Key words} ]
{\small stochastic control, utility maximization, rough volatility, Volterra-Wishart model, Riccati-Volterra equations, non-Markovian}
\end{description}
}
\end{minipage}

\section{Introduction}\label{intro}

Since the observation was made that the paths of realized volatilities are rougher than established volatility models would suggest, cf. \cite{GR18}, there is a growing research interest in developing new models that better fit empirical data. In \cite{EER19}, the popular Heston model \cite{Hest93} was adapted to the rough volatility framework by using a fractional process with Hurst index $H<\frac{1}{2}$ as driver of the volatility process. A more general class of volatility models covering the rough Heston model in \cite{EER19} is obtained by modelling the volatility process as a stochastic Volterra equation of convolution type \cite{AJ19, KLP18,PJM21}. Although most of the literature about rough volatility is concerned with option pricing, there are some recent works dealing with Merton portfolio optimization in such models. While \cite{HW20} and \cite{PJM21} are dealing with the Markowitz portfolio problem, the Merton portfolio problem is studied in \cite{FH19,BD20,HW21}. 

Merton's portfolio problem aims at maximizing an investor's utility from terminal wealth with respect to his utility function. The problem for the classical Heston stochastic volatility model was explicitly solved in \cite{Kra05}, based on the represenation result of \cite{Zar01}, and solutions for affine stochastic volatility models were derived in \cite{KMK10}. In \cite{BL13}, the Merton problem was studied for a multi-asset financial market where the volatility is modeled by a matrix-valued Wishart process, using stochastic control theory. 
In the rough framework it is no longer possible to apply the classical stochastic control approach deriving the corresponding Hamilton-Jacobi-Bellman equation, due to the non-Markovianity of the rough volatility processes. In order to circumvent this problem, in \cite{BD20}, Bäuerle and Desmettre use a finite dimensional approximation of the volatility process in order to cast the problem back into the classical framework. However, this only yields explicit solutions in case that there is no correlation between stock and volatility. Han and Wong \cite{HW21} overcome this difficulty using a martingale distortion Ansatz and applying the martingale optimality principle to obtain an explicit solution for the optimal investment strategy in a mono-asset Volterra Heston model.
In order to take into account several important stylized facts about real financial markets such as choice among multiple assets, roughness of the volatility, correlation between different stocks and leverage effects, i.e. correlation between a stock and its volatility, multivariate rough volatility models have recently been developed  (cf. \cite{CT19}, \cite{RT19}, \cite{AJ19}).
In \cite{PJM21}, Abi Jaber et al. study the Markowitz portfolio problem for a class of multivariate affine Volterra models, that features correlation between the stocks and between a stock and its volatility. 

In this paper we solve the Merton portfolio problem for an investor with a power-utility function for different multivariate Volterra models including the rough Heston model. 
The outline of the paper is as follows: 
Section \ref{Calculus} gives an overview of the calculus of convolutions and resolvents which is needed throughout the paper.
In Section \ref{model} we introduce a class of multivariate affine Volterra models studied in \cite{AJ19} and \cite{PJM21}. For such a market model we consider two different approaches to solve the Merton portfolio problem. We first adapt the martingale distortion transformation used in \cite{HW21} to the multivariate case. However, as it is pointed out in \cite{PJM21}, this only works if the correlation structure is highly degenerate.
Inspired by the techniques used in \cite{BL13}, we then provide a solution for the Merton portfolio problem for a more general correlation structure using calculus of convolutions and resolvents.
In {Section \ref{Wishart}} we introduce a more general market model where the volatility is a matrix-valued stochastic process. In our model we adapt the Wishart stochastic volatility model studied in \cite{GS10} and \cite{BL13} to the Volterra framework, defining the variance-covariance matrix as the solution of a matrix-valued Volterra-Wishart equation, thus extending the Heston model. Considering a matrix-valued volatility process allows us to take into account correlation between different assets.
Despite the non-Markovianity of these settings, the optimal strategy can be expressed explicitly in terms of the solution of a multivariate Riccati-Volterra equation. In Section \ref{sec: example} we illustrate our results with a numerical example. Section \ref{Existence} provides existence and uniqueness results for the appearing Volterra equations even in the matrix-valued case.  
Auxiliary results and longer proofs can be found in the Appendix.

\section{Convolutions and resolvents}\label{Calculus}

In this section we give a short overview of some important definitions and results from the calculus of convolutions and resolvents, that we are going to use frequently throughout the paper. We start by defining three different types of convolutions:

\begin{de}[Convolution of two functions]\cite[Chapter 2]{AJ19}
Let $K$ and $F$ be functions defined on $\R_+$. Then the convolution $K\ast F$ of $K$ and $F$ is defined as
\begin{equation}
(K\ast F)(t)=\int_0^t K(t-s)F(s)ds,
\end{equation}
whenever the above expressions are well-defined.
\end{de}
This definition can of course be extended to matrix-valued functions. In this case it is important that the dimensions of the matrices are compatible.

\begin{de}[Convolution of a measurable function and a measure]\cite[Chapter 2]{AJ19}
Let $K$ be a measurable function on $\R_+$ and $L$ be a measure on $\R_+$ of locally bounded variation. Then the convolutions $K\ast L$ and $L\ast K$ are defined as
\begin{equation}
(K\ast L)(t)=\int_{[0,t]}K(t-s)L(ds);\quad (L\ast K)(t)=\int_{[0,t]}L(ds)K(t-s),
\end{equation}
for all $t\in \R_+$ for which the above integrals exist.
\end{de}

\begin{de}[Convolution of a measurable function and a local martingale]\cite[Chapter 2]{AJ19}
Let $M$ be a $d$-dimensional continuous local martingale and $K:\R_+\rightarrow \R^{m\times d}$ be a function. Then the convolution $K\ast dM$ is defined as
\begin{equation}
(K\ast dM)_t=\int_0^t K(t-s)dM_s.
\end{equation}
\end{de}

\begin{rem}
The above convolution is well-defined as an Itô-integral for any $t\geq 0$ satisfying
\[
\int_0^t \abs{K(t-s)}^2 d\operatorname{tr}(M_s) <\infty.
\]
If $K\in L_{\operatorname{loc}}^2(\R_+)$ and there is a locally bounded process $a$ such that $\langle M \rangle_s=\int_0^s a_u du$ then the convolution is well defined for all $t\geq 0$.
\end{rem}

The following lemma shows that also this type of convolution is associative.

\begin{lem}\cite[Lemma 2.1]{AJ19}\label{lem:assoc}
Let $K\in L_{\operatorname{loc}}^2(\R_+,\R^{m\times d})$ and let $L$ be an $\R^{n\times m}$-valued measure on $\R_+$ of locally bounded variation. Let $M$ be a $d$-dimensional continuous martingale with $\langle M\rangle_t=\int_0^t a_s ds $, $t\geq 0$, for some locally bounded adapted process $a$. Then
\begin{equation}\label{eq assoz}
    (L\ast(K\ast dM))_t=((L\ast K)\ast dM)_t
\end{equation}
for every $t\geq 0$. In particular, taking $F\in L_{\operatorname{loc}}^1(\R_+)$ we may apply \eqref{eq assoz} with $L(dt)=F dt$ to obtain $(F\ast(K\ast dM))_t=((F\ast K)\ast dM)_t$.
\end{lem}

Another useful concept related to the integral kernel $K$ are so called resolvents. We distinguish between resolvents of the first and the second kind.

\begin{de}
Let $K\in L_{\operatorname{loc}}^1(\R_+,\R^{d\times d})$ and $L$ be a $\R^{d\times d}$-valued measure on $\R_+$. Then $L$ is called the \emph{resolvent of the first kind} of $K$ if  
\begin{equation}
K\ast L =L\ast K \equiv I
\end{equation}
where $I$ is the $d$-dimensional identity matrix.
\end{de}

\begin{de}
For a kernel $K\in L_{\operatorname{loc}}^1(\R_+,\R^{d\times d})$, $R\in L_{\operatorname{loc}}^1(\R_+,\R^{d\times d})$ is called the \emph{resolvent of the second kind} of $K$ if 
\begin{equation}
K\ast R =R\ast K = K-R.
\end{equation}
\end{de}

The following table (cf. \cite{AJ19}) gives an overview of some kernels and their corresponding resolvents of the first and second kind.

\begin{table}[h!]
\centering
 \begin{tabular}{||c c c c||} 
 \hline
 Type & $K(t)$ & $R(t)$ & $L(dt)$ \\ [0.5ex] 
 \hline\hline
 \text{Constant} & $c$ & $ce^{-ct}$ & $c^{-1}\delta_0(dt)$ \\ 
 \text{Fractional} & $c\frac{t^{\alpha-1}}{\Gamma(\alpha)}$ & $ct^{\alpha-1}E_{\alpha,\alpha}(-ct^{\alpha})$ & $c^{-1}\frac{t^{-\alpha}}{\Gamma(1-\alpha)}dt$ \\
 \text{Exponential} & $ce^{-\lambda t}$ & $ce^{-\lambda t}e^{-ct}$ & $c^{-1}(\delta_0(dt)+\lambda dt)$ \\
 \text{Gamma} & $ce^{-\lambda t}\frac{t^{\alpha-1}}{\Gamma(\alpha)}$ & $ce^{-\lambda t}t^{\alpha-1}E_{\alpha,\alpha}(-ct^{\alpha})$ & $c^{-1}\frac{1}{\Gamma(1-\alpha)}e^{-\alpha t}\frac{d}{dt}(t^{-\alpha}\ast e^{\lambda t})(t)dt$ \\ [1ex] 
 \hline
 \end{tabular}
\end{table}

For a more detailed discussion of the topic we refer to \cite{GLS90}.

\section{A class of multivariate affine Volterra models}\label{model}

To start our investigation, we use the affine Volterra model introduced in \cite[Chapter 4]{PJM21}.
Let $K=\diag(K_1,\dots,K_d)$ with scalar kernels $K_i\in L^2([0,T],\R)$ on the diagonal. 
In our model we consider $d$ stocks and we assume that the price of the $i$th stock has dynamics 

\begin{equation}\label{stock}
dS_t^{i}=S_t^i(r_t+\theta_i V_t^i)dt+S_t^i\sqrt{V_t^i}dW_{1t}^i,
\end{equation}
where $W_{1t}$ is a $d$-dimensional Brownian motion and $\theta_i\geq 0$.
For $N=\diag(\nu_1,\dots, \nu_d)$ and $D\in \R^{d\times d}$ such that $D_{ij}\geq 0$ if $i\neq j$, the volatility $V=(V^1,\dots,V^d)^{\top}$ is defined as a Volterra square-root process

\begin{equation}\label{Vola}
V_t=v_0(t)+\int_0^t K(t-s)DV_s ds+\int_0^t K(t-s)N{\sqrt{\diag(V_s)}}dB_s.
\end{equation}
Here $v_0:\R_+\rightarrow \R_+^d$ is a deterministic function and $B$ is a $d$-dimensional Brownian motion for which the correlation structure with $W_{1}$ is given by
\begin{equation}\label{correlation structure}
dB_s^i=\rho_i dW_{1s}^i + \sqrt{1-\rho_i^2}dW_{2s}^i,\, i=1,\dots,d,
\end{equation}
where $W_{2}$ is a $d$-dimensional Brownian motion independent of $W_{1}$ and $(\rho_1,\dots,\rho_d)\in[-1,1]^d$. In accordance with \cite{PJM21}, we assume that there exists a continuous $\R_{+}^{2d}$-valued weak solution $(V,S)$ to \eqref{stock}-\eqref{Vola} on some filtered probability space $(\Omega,\mathcal{F}, (\mathcal{F})_{t\geq 0},\mathbb{P})$, satisfying the usual conditions.
A function $f$ is completely monotone on $(0, \infty)$ if it is infinitely differentiable on $(0,\infty)$ and $(-1)^n f^n(t)\geq 0$ for all $n\geq 1$ and $t> 0$.
Under the assumption that for each $i=1,\dots,d$, $K_i$ is completely monotone on $(0,\infty)$ and that there exists $\kappa_i\in (0,2]$ and $k_i>0$ such that
\begin{equation}\label{kernel}
\int_0^h K_i^2(t)dt+\int_0^T(K_i(t+h)-K_i(t))^2 dt\leq k_ih^{\kappa_i},\, h>0,
\end{equation}
the existence of a unique in law $\R_+^d$-valued continuous weak solution $V$ of equation \eqref{Vola} is ensured by \cite[Theorem 6.1]{AJ19} in case that $v_0(t)=V_0+\int_0^t K(t-s)b^0 ds$ for some $V_0,b^0\in\R_+^d$ (cf. \cite[Remark 4.1]{PJM21}). For a discussion about existence of a solution for more general input curves $v_0(t)$, see \cite{AJEE}. Note that condition \eqref{kernel} is fulfilled for constant, non-negative kernels, fractional kernels of the form $\frac{t^{H-\frac{1}{2}}}{\Gamma(H+\frac{1}{2})}$ with $H\in (0,\frac{1}{2}]$, and exponentially decaying kernels $e^{-\beta t}$ with $\beta>0$. The existence of $S$ defined via equation \eqref{stock} follows from that of $V$.

\subsection{The optimization problem}

A portfolio strategy $\pi_t=(\pi_{t,1},\dots, \pi_{t,d})$ is an $(\R^d)^*$ valued, progressively measurable process, where $\pi_{t,k}$ represents the proportion of wealth invested into stock $k$ at time $t$. Under a fixed portfolio strategy, the wealth process $(X_t^{\pi})$ has dynamic

\begin{equation}\label{wealth}
dX_t^{\pi}=X_t^{\pi}(r_t+\pi_t{\diag(V_t)}{\theta}^{\top})dt+X_t\pi_t\sqrt{\diag(V_t)}dW_{1t},
\end{equation}
where ${\theta}=(\theta_1,\dots,\theta_d)$.
By $\mathcal{A}$ we denote the set of admissible portfolio strategies. The conditions under which we consider a strategy to be admissible will be specified later. 
We want to solve the Merton portfolio optimization problem for power utility, i.e. our aim is to find the value function $\mathcal{V}(x_0,v_0)$ such that
\begin{equation}\label{problem}
\mathcal{V}(x_0,v_0)=\sup_{\pi\in\mathcal{A}}\EX_{x_0,v_0}[\frac{1}{\gamma}(X_T^{\pi})^{\gamma}];\quad 0<\gamma<1,
\end{equation}
where $\EX_{x_0,v_0}$ is the conditional expectation given $X_0=x_0, V_0=v_0$. The parameter $\gamma$ represents the relative
risk aversion of the investor. Smaller $\gamma$ correspond to higher risk aversion. 
A portfolio strategy $\pi^*$ for which the supremum is attained is called an optimal strategy. Seen as an optimization problem with state process $(X_t)$ this problem is non-Markovian and the standard stochastic control approach cannot be applied.

\subsection{The martingale distortion transformation}

Consider a one-dimensional market model where the risky asset $S_t$ is given by
\[
dS_t=(r_t+\mu(Y_t))S_tdt+\sigma(Y_t)S_tdW_t
\]
and $Y_t$ is a Markovian process defined via the SDE
\[
dY_t=k(Y_t)dt+h(Y_t)dW_t^Y,
\]
where $W$ and $W^Y$ have correlation $\rho$. In this setup, a candidate for the value function
\[
\mathcal{V}(t,x,y):=\sup_{\pi}\EX[\frac{(X_T^{\pi})^{\gamma}}{\gamma}|X_t=x,Y_t=y]
\]
can be obtained by solving the corresponding Hamilton-Jacobi-Bellman equation. The distortion transformation introduced in \cite{Zar01} uses the Ansatz 
\[
\mathcal{V}(t,x,y)=\frac{x^{\gamma}}{\gamma}\Phi(t,y)^c,
\]
where the constant $c$ is defined as $c:=\frac{1-\gamma}{1-\gamma+\gamma\rho^2}$.
With this choice of $c$, the quadratic terms $(\Phi_y)^2$ in the HJB equation cancel out, leading to a linear PDE for $\Phi$
\[
\Phi_t+\Big( \frac{1}{2}h^2\partial_{yy} + k(y)\partial_y +\frac{\gamma}{1-\gamma}\lambda(y)ch(y)\partial_y \Big)\Phi + \frac{\gamma}{c} \Big(r_t+\frac{\lambda^2(y)}{2(1-\gamma)}\Big)\Phi=0,\, \Phi(T,y)=1,
\]
where the Sharpe ratio $\lambda$ is defined as $\lambda(y):=\mu(y)/\sigma(y)$. By the Feynman-Kac Theorem, $\Phi$ can be written as
\[
\Phi(t,y)=\tilde{\EX}[\operatorname{exp}\Big\{\int_t^T \frac{\gamma}{c} \Big(r_s+\frac{\lambda^2(Y_s)}{2(1-\gamma)}\Big)ds\Big\}|Y_t=y],
\]
where under the probability measure $\tilde{\mathbb{P}}$ with Radon-Nikodym density 
 \[\frac{d\tilde{\mathbb{P}}}{d\mathbb{P}}|_{\mathcal{F}_t}=\operatorname{exp}\Big(\int_0^t \frac{c\gamma}{1-\gamma}\lambda(Y_s)dW_s-\frac{1}{2}\int_0^t\frac{c^2\gamma^2}{(1-\gamma)^2}\lambda^2(Y_s)ds\Big),
 \]
 $\tilde{W}_t^Y=W_t^Y-\int_0^t \frac{c\gamma}{1-\gamma}\lambda(Y_s)ds$ is a standard Brownian motion.

In \cite{FH19} Fouque and Hu showed that if the Sharpe-ratio $\lambda$ is bounded and has bounded derivative, then the value process $\mathcal{V}_t$ can be expressed as $\mathcal{V}_t(x,y)=J_t(X_t=x,Y_t=y)$, where 
\[
J_t(X_t^{\pi},Y_t):=\frac{(X_t^{\pi})^{\gamma}}{\gamma}\Big(\tilde{\EX}\Big[\operatorname{exp}\Big\{\int_t^T \frac{\gamma}{c} \Big(r_s+\frac{\lambda^2(Y_s)}{2(1-\gamma)}\Big)ds\Big\}|\mathcal{F}_t\Big]\Big)^c
\]
even if the volatility process $Y_t$ is non-Markovian.
This approach is called the martingale distortion transformation and was first introduced in the seminal paper \cite{Zar01} and later transferred to a non-Markovian setting in \cite{Teh04}. The extension to the multi-asset case is straight forward in the case of a bounded risk premium (cf. \cite{FH19}, Remark 2.5.).

\subsection{The degenerate correlation case}\label{degenerate}
In this section we present an extension of the proof of \cite{HW21} to the multivariate case for a degenerate correlation structure, i.e. we assume that the correlation in \eqref{correlation structure} is of the form $(\rho,\dots, \rho)$ for $\rho\in [-1,1]$. Note that since in our model the risk premium is unbounded, we can not apply the results of \cite{FH19}. As in the one dimensional case, the Ansatz 
\[
J_t^{\pi}=\frac{(X_t^{\pi})^{\gamma}}{\gamma}\Big(\tilde{\EX}\Big[\operatorname{exp}\Big\{\int_t^T \frac{\gamma}{c} \Big(r_s+\frac{\theta \diag(V_s) \theta^{\top}}{2(1-\gamma)}\Big)ds\Big\}|\mathcal{F}_t\Big]\Big)^c.
\]
is inspired by the martingale distortion transformation described in the previous section. 
Here we use the short notation $J_t^{\pi}$ for $J_t(X_t^{\pi},Y_t)$. Define the diagonal matrices $P:=\diag{(\rho_1,\dots,\rho_d)}$, $\Theta:=\diag{(\theta_1,\dots,\theta_d)}$, $\Psi:=\diag{(\psi_1,\dots,\psi_d)}$ and recall that $N=\diag(\nu_1,\dots, \nu_d)$.
Under the new probability measure $\tilde{\mathbb{P}}$ defined via the Radon-Nikodym density 
\[
\frac{d\tilde{\mathbb{P}}}{d\mathbb{P}}|_{\mathcal{F}_t}=\operatorname{exp}\Big(\frac{\gamma}{1-\gamma}\int_0^t \theta \sqrt{\operatorname{diag}({V_s}})dW_{1s}-\frac{\gamma^2}{2(1-\gamma)^2}\int_0^t \theta {\operatorname{diag(V_s)}}\theta^{\top}\Big)
\]
together with the new standard brownian motion under $\tilde{\mathbb{P}}$ 
\[
\tilde{W}_{1t}=W_{1t}-\frac{\gamma}{1-\gamma}\int_0^t \theta \sqrt{\operatorname{diag}(V_s)}ds,
\]
 an application of the exponential-affine transform formula in \cite[Theorem 4.3]{AJ19} yields
\begin{equation}\label{M}
    \begin{aligned}
 {}
    & \tilde{\EX}\Big[\operatorname{exp}\Big\{\int_t^T \Big(\frac{\gamma\theta\Theta V_s }{2c(1-\gamma)}\Big)ds\Big\}|\mathcal{F}_t\Big]=\\
    &\qquad=\operatorname{exp}\Big\{\int_t^T \frac{\gamma \theta\Theta\xi_t(s)}{2(1-\gamma)}+\frac{c}{2}\psi(T-s)N^2\Psi(T-s)\xi_t(s)ds \Big\}=:M_t,
    \end{aligned}
\end{equation}
where $\xi_t(s):={\tilde{\EX}[V_s|\mathcal{F}_t]}$ denotes the conditional $\tilde{\mathbb{P}}$-expected variance and $\psi\in L^2([0,T],(\R^d)^*)$ solves the Ricatti-Volterra equation
\[
\psi= \big(\frac{\gamma}{2c(1-\gamma)}\theta\Theta+\psi\Lambda +\frac{1}{2}\psi N^2\Psi\big)\ast K,
\]
with $\Lambda=D+\frac{\gamma}{1-\gamma}NP\Theta$.
Thus we obtain
\[
J_t^{\pi}=\frac{(X_t^{\pi})^{\gamma}}{\gamma}M_t.
\]

In order to find the value function and the optimal strategy, we show that the family $\{J_t^{\pi}\}_{\pi\in\mathcal{A}}$ fulfills the \emph{martingale optimality principle} (cf. \cite{HI05,JS12,HW21}), i.e. we show that:
\begin{enumerate}
    \item $J_T^{\pi}=\frac{1}{\gamma}(X_T^{\pi})^{\gamma}$ for all $\pi\in\mathcal{A}$;
    \item $J_0^{\pi}=J_0$ is a constant independent of $\pi$;
    \item $J_t^{\pi}$ is a supermartingale for all $\pi\in\mathcal{A}$ and there exists $\pi^*\in\mathcal{A}$ such that $J_t^{\pi^*}$ is a martingale.
\end{enumerate}
A family of processes with the above properties can now be used to compare the expected utilities of an arbitrary strategy $\pi$ and the strategy $\pi^*$:
\begin{equation*}
    \EX[\frac{1}{\gamma}(X_T^{\pi})^{\gamma}]=\EX[J_T^{\pi}]\leq J_0^{\pi}=J_0^{\pi^*}=\EX[J_T^{\pi^*}]=\EX[\frac{1}{\gamma}(X_T^{\pi^*})^{\gamma}]=\mathcal{V}(x_0,v_0).
\end{equation*}
Thus the strategy $\pi^*$ is indeed our desired optimal portfolio strategy.

\begin{de}\label{Def:adm}
In the setting described above, we say that a portfolio strategy $\pi$ is admissible if
\begin{enumerate}
    \item the SDE \eqref{wealth} for the wealth process $(X_t^{\pi})$ has a unique solution in terms of $(S,V,W_1)$;
    \item $\EX[\frac{1}{\gamma}(X_T^{\pi})^{\gamma}]<\infty$ for all $0<\gamma<1$;
    \item $\int_0^t \pi_s\diag{(V_s)}\pi_s^{\top}ds<\infty$ a.s.
\end{enumerate}

\end{de}

The main result we get for the degenerate correlation case is the following, our proof enhances the arguments of the proof of \cite[Theorem 3.3]{HW21} to the multi-dimensional case:

\begin{thm}
Let $\Lambda=D+\frac{\gamma}{1-\gamma}NP\Theta$ be invertible and let $\psi$ be the unique, continuous non-continuable solution of the Riccati-Volterra equation
\begin{equation}\label{RV1}
\psi(t)=\int_0^t F_1(\psi)(t-s)K(s)ds
\end{equation}
\begin{equation}\label{RV3}
F_1(\psi)= \frac{\gamma}{2c(1-\gamma)}\theta\Theta+\psi\Lambda +\frac{1}{2}\psi N^2\Psi
\end{equation}
on the interval $[0,T_{\operatorname{max}}]$, given by Theorem \ref{vecRicEx}\footnote{{More details on $T_{max}$ can be found in Section 6.2.}}. Then $J_t^{\pi}=\frac{(X_t^{\pi})^{\gamma}}{\gamma}M_t$ satisfies the martingale optimality principle for $t\in [0,T]$, $T\leq T_{\operatorname{max}}$ and the optimal portfolio strategy $\pi^*$ is given by
\begin{equation}
    \pi_t^*=\frac{1}{1-\gamma}(\theta+c\psi(T-t)NP).
\end{equation}
\end{thm}

\emph{Proof:}
We show that $J_t^{\pi}$ fulfills the martingale optimality principle. 
For the first condition, note that $M_T=1$ and hence $J_T^{\pi}=\frac{1}{\gamma}(X_T^{\pi})^{\gamma}$. Since $M_0$ is a constant independent of $\pi$, $J_0^{\pi}=\frac{x_0^{\gamma}}{\gamma}M_0$ is also independent of $\pi$ and thus also the second condition is satisfied.
In order to show that also the third condition is fulfilled, we apply Itô's formula on $J_t^{\pi}$. Using Lemma \ref{dynamics M}, this yields

\begin{equation*}
    \begin{aligned}
    dJ_t^{\pi}{}
    &= (r_t+\pi_t\diag{(V_t)}\theta^{\top})M_tX_t^{\gamma}dt + M_tX_t^{\gamma}\pi_t\sqrt{\diag(V_t)}dW_{1t}\\
    &\quad - M_tX_t^{\gamma}(r_t+\frac{1}{2(1-\gamma)}\theta\diag(V_t)\theta^{\top})dt - \frac{1}{1-\gamma}M_tX_t^{\gamma}c\psi(T-t)NP_1\diag(V_t)\theta^{\top}dt\\
    &\quad - \frac{1}{2(1-\gamma)}M_tX_t^{\gamma}c^2\psi(T-t)N^2P_1^2\Theta^2\diag{(V_t)}\psi^{\top}(T-t)dt \\
    &\quad + M_t\frac{X_t^{\gamma}}{\gamma}c\psi(T-t)NP_1\sqrt{\diag(V_t)}dW_{1t} + M_t\frac{X_t^{\gamma}}{\gamma}c\psi(T-t)NP_2\sqrt{\diag(V_t)}dW_{2t}\\
    &\quad + \frac{\gamma-1}{2}M_tX_t^{\gamma}\pi_t\diag{(V_t)}\pi_t^{\top}dt + M_tX_t^{\gamma}c\pi_t NP_1\diag(V_t)\psi^{\top}(T-t)dt\\
    &= J_t^{\pi}F(\pi,t)dt+J_t^{\pi}(c\psi(T-t)NP_1\sqrt{\diag(V_t)}+\gamma\pi_t\sqrt{\diag{(V_t)}})dW_{1t}\\
    & \quad + J_t^{\pi}c\psi(T-t)NP_2\sqrt{\diag(V_t)}dW_{2t}
    \end{aligned}
\end{equation*}

with 

\begin{equation*}
    \begin{aligned}
    F(\pi,t){}
    &= \frac{\gamma(\gamma-1)}{2}\pi_t\diag(V_t)\pi_t^{\top}+\gamma\pi_t(\diag(V_t)\theta^{\top}+cNP\diag{(V_t)}\psi^{\top}(T-t))\\
    & - \frac{\gamma}{2(1-\gamma)}\lVert \diag(V_t)\theta^{\top}+cNP_1\diag{(V_t)}\psi^{\top}(T-t) \rVert_2^2. 
    \end{aligned}
\end{equation*}

Note that $F(\pi^{*},t)=0$ and since $F$ is a quadratic function in $\pi$ and $\gamma\in(0,1)$ we have $F(\pi,t)\leq 0$. 
Solving the stochastic differential equation for $J_t^{\pi}$ yields
\[
J_t^{\pi}= \frac{M_0 x_0^{\gamma}}{\gamma}e^{\int_0^t F(\pi_s,s)ds} G(\pi_t,t)
\]
with

\begin{equation*}
    \begin{aligned}
    G(\pi,t){}
    &= \operatorname{exp}\{-\frac{1}{2}\int_0^t (\lVert cNP_1\diag{(V_t)}\psi^{\top}(T-t) + \gamma \sqrt{\diag(V_t)}\pi_t^{\top} \rVert_2^2 \\
    &\qquad\qquad\qquad + \lVert cNP_2\diag{(V_t)}\psi^{\top}(T-t)  \rVert_2^2) ds\\
    &\quad+ \int_0^t[c\psi(T-s)NP_1\sqrt{\diag(V_s)}+\gamma \pi_s\sqrt{\diag(V_s)}]dW_{1s}\\
    &\quad+ \int_0^t c\psi(T-s)NP_2\sqrt{\diag(V_s)}dW_{2s}\}.
    \end{aligned}
\end{equation*}

Now, since $F(\pi,t)\leq 0$, $e^{\int_0^t F(\pi_s,s)ds}$ is a non-increasing function. By our assumptions on the admissible strategies, $\int_0^t \pi_s\diag{(V_s)}\pi_s^{\top}ds<\infty$ and thus the stochastic exponential $G$ is a local martingale (which follows from the basic properties of the Dool\'{e}ans-Dade exponential). Therefore we can find a sequence of stopping times $\tau_1,\tau_2,\tau_3,\dots$ with $\operatorname{lim}_{n\rightarrow\infty} \tau_n=T$ a.s. satisfying 
\[
\EX[J_{t\wedge \tau_n}^{\pi}|\mathcal{F}_s]\leq J_{s\wedge \tau_n}^{\pi},\, s\leq t.
\]
Since $J_t^{\pi}\geq 0$, an application of Fatou's Lemma yields that $J_t^{\pi}$ is a supermartingale for every arbitrary admissible strategy $\pi$.
It remains to show that $J_t^{\pi^*}$ is a true martingale for the optimal strategy $\pi^*$. In this case $e^{\int_0^t F(\pi_s,s)ds}=1$ and hence $J_t^{\pi^*}=\frac{M_0 x_0^{\gamma}}{\gamma} G(\pi_t^*,t)$. $G(\pi^*,t)$ is a martingale by Lemma \ref{extended AJ lemma} and so is $J_t^{\pi^*}$.
In order to show that $\pi^*$ is admissible, we have to show that (a), (b), (c) of Definition~\ref{Def:adm} hold. Part (a) is true because \eqref{wealth} has a unique solution in terms of $(S,V,W_1)$ as $\pi^*$ is deterministic.
For part (b) it suffices to show that
\[\mathbb{E}[\operatorname{sup}_{t\in[0,T]}|X_t^{\pi^*}|^q]<\infty.
\]
Inserting the explicit solution of the wealth process $X^{\pi^*}$ into the left-hand side and applying a combination of Doob's maximal inequality together with Hölder's inequality then leads to the desired result. Part (c) is true as our resulting optimal strategies are deterministic.\qed

\subsection{The general correlation case}\label{general}

For the case where the correlation in \eqref{correlation structure} is given by an arbitrary vector $(\rho_1,\dots,\rho_d)\in[-1,1]^d$, the martingale distortion arguments from the previous section do not work anymore. Therefore, we develop a new approach inspired by the verification arguments used in \cite{BL13} to solve the optimization problem for this more general correlation structure. In this setting we say that a portfolio strategy $\pi$ is admissible if
\begin{enumerate}
    \item the SDE \eqref{wealth} for the wealth process $(X_t^{\pi})$ has a unique strong solution;
    \item $\EX[\frac{1}{\gamma}(X_T^{\pi})^{\gamma}]<\infty$ for all $0<\gamma<1$;
    \item $\pi$ is bounded.
\end{enumerate}
Remember that $N=\diag{(\nu_1,\dots,\nu_d)}$, $P=\diag{(\rho_1,\dots,\rho_d)}$, $\Theta=\diag{(\theta_1,\dots,\theta_d)}$, $\Psi=\diag{(\psi_1,\dots,\psi_d)}$.
The main result we provide for this case is the following:

\begin{thm}
For $\Lambda=D+\frac{\gamma}{1-\gamma}NP\Theta$, let $\psi$ be the solution of the Riccati-Volterra equation
\begin{equation}\label{RV2}
\psi(t)=\int_0^t F_2(\psi)(t-s)K(s)ds
\end{equation}
with
\begin{equation}\label{RV4}
F_2(\psi)=\frac{\gamma}{2(1-\gamma)}{\theta}\Theta+\psi\Lambda +\frac{1}{2}(\psi N^2\Psi+\frac{\gamma}{1-\gamma}\psi N^2P^2\Psi)
\end{equation}
on the interval $[0,T_{\operatorname{max}})$, given by Theorem \ref{vecRicEx}. Then for $t\in [0,T]$, $T< T_{\operatorname{max}}$, an optimal investment strategy $\pi_t^*$ for the Merton portfolio problem \eqref{problem} is given by 
\begin{equation}\label{optimal strategy}
    \pi_t^*=\frac{1}{1-\gamma}(\theta+\psi(T-t)NP)
\end{equation}
and the value function can be written as
\[
\mathcal{V}(x_0,V_0)=\frac{x_0^\gamma}{\gamma}\expo\Big(\int_0^T \gamma r_s +  F_2(\phi)(T-s)v_0(s)ds\Big).
\]
\end{thm}

\textit{Proof:}
The proof is a straight forward adaption of the arguments from the proof of Theorem \ref{mainthm} to the vector-valued case. \qed

\subsection{Comparison of the different approaches}\label{comparison}
The martingale distortion approach in Section \ref{degenerate} yields the following solution for the Merton portfolio problem in the $d$-dimensional affine Volterra model.
The optimal portfolio strategy is given by
\[
\pi_t^*=\frac{1}{1-\gamma}(\theta+c\psi(T-t)NP).
\]
The value function can be written as
\[
\sup_{\pi\in \mathcal{A}_1}\EX^{t,X_t,V_t}[\frac{1}{\gamma}(X_T^{\pi})^\gamma] = H(t,X_t,V_t)
\]
with
\[
H(t,X_t,V_t)=\frac{X_t^\gamma}{\gamma}\expo\Big(\int_t^T(\gamma r_s + \frac{\gamma}{2(1-\gamma)}\theta\Theta\xi_t(s)+\frac{c}{2}\psi(T-s)N^2\Psi(T-s)\xi_t(s))ds\Big)
\]
where $\psi$ is the solution of the Riccati-Volterra equation
\begin{equation}
\psi(t)=\int_0^t F_1(\psi)(t-s)K(s)ds
\end{equation}
\[
F_1(\psi)= \frac{\gamma}{2c(1-\gamma)}\theta\Theta+\psi\Lambda +\frac{1}{2}\psi N^2\Psi.
\]
For the approach in Section \ref{general} we get the following results.
The optimal portfolio strategy is given by
\[
\pi_t^*=\frac{1}{1-\gamma}(\theta+\phi(T-t)NP).
\]
The value function can be written as
\[
\sup_{\pi\in \mathcal{A}_2}\EX^{0,x_0,V_0}[\frac{1}{\gamma}(X_T^{\pi})^\gamma] = G(0,x_0,V_0)
\]
with
\[
G(0,x_0,V_0)=\frac{x_0^\gamma}{\gamma}\expo\Big(\int_0^T \gamma r_s +  F_2(\phi)(T-s)V_0(s)ds\Big)
\]
where $\phi$ is the (unique, global) solution of the Riccati-Volterra equation
\begin{equation}
\phi(t)=\int_0^t F_2(\phi)(t-s)K(s)ds
\end{equation}
\[
F_2(\phi)=\frac{\gamma}{2(1-\gamma)}{\theta}\Theta+\phi\Lambda +\frac{1}{2}(\phi N^2\Phi+\frac{\gamma}{1-\gamma}\phi N^2P^2\Phi).
\]
\begin{lem}
Let $\rho_1=\dots =\rho_d$ and define $c:=\frac{1-\gamma}{1-\gamma+\gamma\rho^2}$. If $\psi$ is the unique global solution of the Riccati Volterra equation
\[
\psi(t)=\int_0^t F_1(\psi)(t-s)K(s)ds,
\]
then $\phi=c\psi$ is the unique global solution of
\[
\phi(t)=\int_0^t F_2(\phi)(t-s)K(s)ds.
\].
\end{lem}
\emph{Proof.} 
First we show that $cF_1(\psi)=F_2(c\psi)$.
\begin{equation*}
    \begin{aligned}
    cF_1(\psi){}
    &= c(\frac{\gamma}{2c(1-\gamma)}{\theta}\Theta+\psi\Lambda +\frac{1}{2}\psi N^2\Psi) = \frac{\gamma}{2(1-\gamma)}{\theta}\Theta+c\psi\Lambda +\frac{c}{2}\psi N^2\Psi\\
    &= \frac{\gamma}{2(1-\gamma)}{\theta}\Theta+(c\psi)\Lambda +\frac{1}{2c}(c\psi) N^2(c\Psi)\\
    &= \frac{\gamma}{2(1-\gamma)}{\theta}\Theta+(c\psi)\Lambda +\frac{1}{2}\frac{1-\gamma+\gamma\rho^2}{1-\gamma}(c\psi) N^2(c\Psi)\\
    &= \frac{\gamma}{2(1-\gamma)}{\theta}\Theta+(c\psi)\Lambda +\frac{1}{2}[(c\psi) N^2(c\Psi) +\frac{\gamma}{1-\gamma}(c\psi) N^2P^2(c\Psi)]\\
    &= F_2(c\psi)
    \end{aligned}
\end{equation*}
Since $\psi$ is the unique solution of $\psi(t)=\int_0^t F_1(\psi)(t-s)K(s)ds$, we have
\[
c\psi(t)=\int_0^t cF_1(\psi)(t-s)K(s)ds=\int_0^t F_2(c\psi)(t-s)K(s)ds
\]
and thus $\phi=c\psi$ has to be the unique solution of $\phi(t)=\int_0^t F_2(\phi)(t-s)K(s)ds$.
\qed
\begin{thm}
Let $\rho_1=\dots =\rho_d$. Then $\pi^*(\mathcal{A}_1)=\pi^*(\mathcal{A}_2)$ and for the value functions we have $H(0,x_0,V_0)=G(0,x_0,V_0)$.
\end{thm}
\emph{Proof.}
Recall that $\pi_t^*(\mathcal{A}_1)=\frac{1}{1-\gamma}(\theta+c\psi(T-t)NP)$ and $\pi_t^*(\mathcal{A}_2)=\frac{1}{1-\gamma}(\theta+\phi(T-t)NP)$. By the Lemma we have $c\psi=\phi$ and the equality follows immediately.
It remains to show that the value functions are equal. From \cite[Lemma 4.2]{AJ19} we know that 
\[
\xi_0(t)=\EX[V_t]=(I-\int_0^t R_{\Lambda}(u)du)V_0.
\]
Using this fact, the value function $H$ reads
\begin{equation*}
    \begin{aligned}
    {}
    & H(0,x_0,V_0)=\\
    &= \frac{x_0^\gamma}{\gamma}\expo\Big(\int_0^T\gamma r_s + \frac{\gamma}{2(1-\gamma)}\theta\Theta\xi_0(s)+\frac{c}{2}\psi(T-s)N^2\Psi(T-s)\xi_0(s)ds\Big)\\
    &= \frac{x_0^\gamma}{\gamma}\expo\Big(\int_0^T\gamma r_s + cF_1(\psi)(T-s)\xi_0(s)+c\psi(T-s)\Lambda\xi_0(s)ds\Big)\\
    &= \frac{x_0^\gamma}{\gamma}\expo\Big(\int_0^T\gamma r_s + cF_1(\psi)(T-s)(I-\int_0^s R_{\Lambda}(u)du)V_0(s)\\
    &\qquad\qquad\qquad\qquad\qquad\qquad\qquad +c\psi(T-s)\Lambda(I-\int_0^s R_{\Lambda}(u)du)V_0(s)ds\Big)\\
    &= \frac{x_0^\gamma}{\gamma}\expo\Big(\int_0^T\gamma r_s + cF_1(\psi)(T-s)V_0(s)ds\Big)\\
    & \cdot \expo\Big(\int_0^T -cF_1(T-s)\int_0^s R_{\Lambda}(u)du V_0(s) ds + \int_0^T c(F_1(\psi)\ast K)(T-s)\Lambda IV_0(s)ds\\
    &\qquad\qquad\qquad\qquad\qquad\qquad - \int_0^T c(F_1(\psi)\ast K)(T-s)\Lambda \int_0^s R_{\Lambda}(u)du V_0(s)ds\Big)\\
    &= \frac{x_0^\gamma}{\gamma}\expo\Big(\int_0^T\gamma r_s + cF_1(\psi)(T-s)V_0(s)ds\Big) \cdot \expo\Big( -cF_1(\psi)\ast ((R_{\Lambda}\ast I)V_0)(T)\\
    &\qquad\qquad\qquad +c((F_1(\psi)\ast K)\Lambda)\ast IV_0(T)- c((F_1(\psi)\ast K)\Lambda)\ast((R_{\Lambda}\ast I)V_0)(T)\Big)\\ 
    &= \frac{x_0^\gamma}{\gamma}\expo\Big(\int_0^T\gamma r_s + cF_1(\psi)(T-s)V_0(s)ds\Big) \cdot \expo\Big( -(cF_1(\psi)\ast R_{\Lambda}\ast IV_0)(T)\\
    &\qquad\qquad\qquad +(cF_1(\psi)\ast K\Lambda\ast IV_0)(T)- (cF_1(\psi)\ast (K\Lambda\ast R_{\Lambda})\ast IV_0)(T)\Big)\\ 
    &= \frac{x_0^\gamma}{\gamma}\expo\Big(\int_0^T\gamma r_s + cF_1(\psi)(T-s)V_0(s)ds\Big) \cdot \expo\Big( -(cF_1(\psi)\ast R_{\Lambda}\ast IV_0)(T)\\
    &\qquad\qquad\qquad +(cF_1(\psi)\ast K\Lambda\ast IV_0)(T)- (cF_1(\psi)\ast (K\Lambda- R_{\Lambda})\ast IV_0)(T)\Big)\\ 
    &= \frac{x_0^\gamma}{\gamma}\expo\Big(\int_0^T\gamma r_s + cF_1(\psi)(T-s)V_0(s)ds\Big)\cdot \expo(0)\\
    &= \frac{x_0^\gamma}{\gamma}\expo\Big(\int_0^T\gamma r_s + F_2(c\psi)(T-s)V_0(s)ds\Big)=G(0,x_0,V_0).\qed
    \end{aligned}
\end{equation*}

\section{The Volterra-Wishart volatility model}\label{Wishart}
In this section we present the Volterra-Wishart model which is a generalization of the Wishart volatility model described in \cite{GS10} and \cite{BL13} to the Volterra framework. In contrast to the class of models presented in the previous section, the volatility is now modeled as a matrix-valued stochastic process. 
The main advantage of using a matrix-valued volatility model is that this allows us to take into account the correlation between different stocks in our market and therefore extending the vector-valued model presented in the previous section. In contrast to the quadratic volatility models described in \cite[Chapter 5]{PJM21} and \cite{AJ22}, generalizing the Stein-Stein and Schöbel-Zhu model respectively, our model is an extension of the Heston model to the multivariate Volterra framework. The reason why we in particular investigate a Wishart model is that, besides being a straight forward generalisation of the popular Heston model, the additional degrees of freedom with regard to the stochastic correlation enable a better fit of market data while being still efficiently tractable, compare e.g. \cite{BBK08} and \cite{LBM21}. 

\subsection{Market model and optimization problem}

{In our model the market consists of one riskfree asset with time-dependent, deterministic risk free rate $r_t$ and $d$ risky assets.
The asset return vector process $(S_t)_{t\geq 0}=(S_{t,1},\dots, S_{t,d})_{t\geq 0}$ is defined via the stochastic differential system
\begin{equation}\label{assets}
dS_t=\operatorname{diag}(S_t)((\underline{r}_t+\Sigma_t v)dt+\Sigma_t^{1/2}dW_t^S),
\end{equation}
where $(W_t^{S})_{t\geq 0}$ is a $d$-dimensional Brownian Motion vector {, $\underline{r}_t=(r_t,\dots,r_t)$ and $v\in\R^d$}. The stochastic volatility process $(\Sigma_t)_{t\geq 0}$ is given by the solution of the matrix-valued Volterra equation
\begin{equation}\label{vola}
\begin{aligned}
&\Sigma_t=\Sigma_0+\int_0^t K(t-s) (NN^{\top}+M\Sigma_s+\Sigma_sM^{\top})ds\\
&\qquad+\int_0^t K(t-s)\Sigma_s^{1/2}dW_s^{\sigma}Q+\int_0^t K(t-s)Q^{\top}(dW_s^{\sigma})^{\top}\Sigma_s^{1/2}.
\end{aligned}
\end{equation}
The integral kernel in the above equation $K=\diag(K_1,\dots, K_d)$, {$K_1=\dots=K_d$} is diagonal with scalar kernels $K_i\in L^2([0,T],\R)$. The deterministic initial value $\Sigma_0$ is assumed to be positive definite. The $d\times d$ matrices $N$ and $M$ are responsible for the mean reversion, while the matrix $Q$ governs the volatility of the process which is driven by a $d\times d$ Brownian motion matrix $(W_t^{\sigma})_{t\geq 0}$. In order to incorporate the leverage effect in our model, we allow the Brownian motions $(W_t^{S})_{t\geq 0}$ and $(W_t^{\sigma})_{t\geq 0}$ to be correlated and we assume that $d\langle W_{t,k}^S, W_{t,ij}^{\sigma} \rangle=\rho_{k,ij}$ with $\rho_{k,ij}=0$ for $k\neq i$ and $\rho_{k,kj}=:\rho_j$ independent of $k$. Thus for $\rho=(\rho_1,\dots,\rho_d)^{\top}$ and another $d$-dimensional Brownian motion vector $B_t$ independent of $(W_t^{\sigma})$, we have
\[
 W_t^S=\sqrt{1-\rho^{\top}\rho}B_t+W_t^{\sigma}\rho.
\] 
Under the assumption that the components of the kernel $K$ fulfill the condition
  \begin{equation}\label{kernels}
  \begin{aligned}
& K_i\in L_{\operatorname{loc}}^2(\R_+,\R)\text{ and there is } \kappa_i\in(0,2]\text{ such that } \int_0^h K(t)^2 dt=O(h^{\kappa_i})\\
& \text{ and } \int_0^T (K(t+h)-K(t))^2 dt =O(h^{\kappa_i}) \text{ for every } T<\infty,
  \end{aligned}
  \end{equation} 
Theorem \ref{ExVol} ensures the existence of a continuous, symmetric and positive definite $\R^{d\times d}$-valued local weak solution $\Sigma$ to \eqref{assets}-\eqref{vola} on some filtered probability space $(\Omega,\mathcal{F}, (\mathcal{F})_{t\geq 0},\mathbb{P})$, satisfying the usual conditions.
For such a symmetric, positive definite matrix $\Sigma$, by $\Sigma^{1/2}$ we denote the unique symmetric, positive definite matrix $M$ for which $M^2=\Sigma$.}
Note that for $K=I$ we recover the classical Wishart model described in \cite{BL13}. 

A portfolio strategy $\pi_t=(\pi_{t,1},\dots, \pi_{t,d})^{{\top}}$ is an $\R^d$-valued, progressively measurable process, where $\pi_{t,k}$ represents the proportion of wealth invested into stock $k$ at time $t$. Under a fixed portfolio strategy, the wealth process $(X_t^{\pi})$ has dynamics

\begin{equation}\label{wealth1}
dX_t^{\pi}=X_t^{\pi}[(r_t+\pi_t^{\top}\Sigma_t^{}v) dt+\pi_t^{\top}\Sigma_t^{1/2}dW_{t}^S],\quad X_0=x_0.
\end{equation}

By $\mathcal{A}$ we denote the set of admissible portfolio strategies. 
\begin{de}
In our setting we say that a portfolio strategy $\pi$ is admissible if
\begin{enumerate}
    \item the SDE \eqref{wealth1} for the wealth process $(X_t^{\pi})$ has a unique solution {in terms of $(S,\Sigma,W^S)$};
    \item $\EX[\frac{1}{\gamma}(X_T^{\pi})^{\gamma}]<\infty$ for all $0<\gamma<1$;
    \item $\pi$ is bounded { a.s}.
\end{enumerate}
\end{de}

We want to solve the Merton portfolio optimization problem for power utility, i.e. our aim is to find the value function $\mathcal{V}(x_0,\Sigma_0)$ such that
\begin{equation}\label{problem1}
\mathcal{V}(x_0,\Sigma_0)=\sup_{\pi\in\mathcal{A}}\EX_{x_0,\Sigma_0}[\frac{1}{\gamma}(X_T^{\pi})^{\gamma}],\quad 0<\gamma<1,
\end{equation}
where $\EX_{x_0,\Sigma_0}$ is the conditional expectation. Again, the parameter $\gamma$ represents the relative risk aversion of the investor. A portfolio strategy $\pi^*$ for which the supremum is attained is called an optimal strategy. As stated before, seen as an optimization problem with state process $(X_t)$ this problem is non-Markovian and the standard stochastic control approach cannot be applied.

\subsection{The main result}

We solve the Merton portfolio problem for the Volterra-Wishart model using a verification argument inspired by \cite{BL13}. {As it was pointed out in Section \ref{degenerate}, the martingale distortion approach, used by \cite{HW21} in the one dimensional case, can only be applied to the multivariate setting if the correlation structure is highly degenerate, i.e. $\rho_1=\dots = \rho_d$. Since in our model we allow arbitrary correlation vectors $\rho$, we have to use  different techniques. }
{The proof builds on ideas presented in \cite{BL13} for the classical Wishart volatility model, yet we are facing serious technical challenges by doing that. In particular, in our case the stochastic volatility process is of convolution type and hence non-Markovian, and therefore we cannot use Itô's formula in order to show optimality of the candidate for the optimal strategy. In order to overcome this, we have to resort to the calculus of convolutions and resolvents. Note that the techniques we apply in our proof rely on the affine structure of the Volterra-Wishart convolution equation. Since the resulting calculations are rather involved, we present a condensed version of the proof here and the full proof  can be found in Appendix~\ref{A:main_result}.} {Note that in contrast to the vector-valued cased, we cannot prove global existence of a week solution for the Volterra-Wishart process \eqref{vola} and also the Riccati equation \eqref{MRV} can only be solved locally. Therefore we can only solve the Merton problem for a time horizon $T$ within the interval $[0;T_{\operatorname{max}}\wedge T_{\operatorname{pos}}]$.}

\begin{thm}\label{mainthm}
Assume that equation \eqref{vola} has a positive definite, continuous weak solution on the interval $[0,T_{\operatorname{pos}})$. 
Let $\psi$ be the solution of the matrix Riccati-Volterra equation
\begin{equation}\label{MRV}
\psi(t)=\int_0^t f(\psi)(t-s)K(s)ds
\end{equation}
with
\begin{equation}\label{fpsi}
    f(\psi)=\psi\tilde{M}+\tilde{M}^{\top}\psi+2\psi\tilde{Q}^{\top}\tilde{Q}\psi+\tilde{\Gamma},
\end{equation}
\begin{equation*}
\tilde{M}=M+\frac{\gamma}{1-\gamma}Q^{\top}\rho v^{\top},\quad \tilde{Q}^{\top}Q=Q^{\top}Q+\frac{\gamma}{1-\gamma}Q^{\top}\rho\rho^{\top}Q,\quad \tilde{\Gamma}= \frac{\gamma}{2(1-\gamma)}vv^{\top}
\end{equation*}
on the interval $[0,T_{\operatorname{max}})$, given by Theorem \ref{ExRic}. Then for $t\in [0,T]$, $T< T_{\operatorname{max}}\wedge T_{\operatorname{pos}}$, an optimal investment strategy $\pi_t^*$ for the Merton portfolio problem \eqref{problem1} is given by 
\begin{equation}\label{optimal strategy1}
    \pi_t^*=\frac{1}{1-\gamma}(v+2\psi(T-t)Q^{\top}\rho)
\end{equation}
and the value function can be written as
\[
\mathcal{V}(x_0,\Sigma_0)=\frac{x_0^\gamma}{\gamma}\expo\Big(\int_0^T \gamma r_s + \operatorname{Tr}[f(\psi)(T-s)\Sigma_0 +\psi(T-s)NN^{\top}]ds\Big).
\]
\end{thm}
{
\emph{Sketch of Proof:}
In order to prove that $\pi^*$ is indeed the optimal portfolio strategy, we show that for 
\[
G(x_0,\Sigma_0):=\frac{x_0^\gamma}{\gamma}\expo\Big(\int_0^T \gamma r_s + \operatorname{Tr}[f(\psi)(T-s)\Sigma_0 +\psi(T-s)NN^{\top}]ds\Big), 
\]
we have
\begin{enumerate}
    \item $\EX^{x_0,\Sigma_0}{\big[\frac{1}{\gamma}(X_T^{\pi^*})^{\gamma}\big]}=G(x_0,\Sigma_0)$ for $\pi_t^*=\frac{1}{1-\gamma}(v+2\psi(T-t)Q^{\top}\rho)$,\label{a1}
    \item $\EX^{x_0,\Sigma_0}{\big[\frac{1}{\gamma}(X_T^{\pi})^{\gamma}\big]}\leq G(x_0,\Sigma_0)$ for every other admissible strategy.\label{b1}
\end{enumerate}
\textbf{Proof of the equality (a)}:
The SDE for the wealth process can be solved explicitly: 
\[
X_T^{\pi}=X_0^{\pi}\expo\Big(\int_0^T(r_s+\pi_s^{\top}\Sigma_s v-\frac{1}{2}\lVert \pi_s^{\top}\Sigma_s^{1/2} \rVert_{2}^2 )ds + \int_0^T \pi_s^{\top}\Sigma_s^{1/2}dW_s^S\Big);\quad X_0^\pi=x_0.
\]
Introducing a new probability measure ${\Q}$ with Radon-Nikodym density  
\begin{equation}\label{density Q}
{Z}_0:=\frac{d{\Q}}{d\mathbb{P}}|_{\mathcal{F}_0}=\expo\Big(\gamma\int_0^T (\pi_s^{*})^{\top}\Sigma_s^{1/2}dW_{s}^S-\frac{\gamma^2}{2}\int_0^T \lVert (\pi_s^{*})^{\top}\Sigma_s^{1/2} \rVert_{2}^2 ds\Big),
\end{equation}
we obtain
\begin{equation*}
\begin{aligned}
x_0^{-\gamma}\EX_{x_0,\Sigma_0}{(X_T^{\pi^*})^{\gamma}}{}
    & =\EX_{x_0,\Sigma_0}^{\Q}\Big[\expo\Big(\int_0^T \gamma r_s ds+\int_0^T \operatorname{Tr}\big[F_s\Sigma_s\big] ds\Big)\Big],\\
\end{aligned}
\end{equation*}
where the matrix-valued, deterministic process $F_s$ is given by 
\[
F_s:=\gamma v(\pi_s^{*})^{\top}
    +\frac{\gamma(\gamma-1)}{2} \pi_s^*(\pi_s^{*})^{\top}.
\]
Inserting the candidate for the optimal strategy $\pi_t^*=\frac{1}{1-\gamma}(v+2\psi(T-t)Q^{\top}\rho)$, $F_s$ can be written in terms of $f(\psi)$
\begin{equation*}
    \begin{aligned}
    F_s &= f(\psi)(T-s)-\psi(T-s)M-M^{\top}\psi(T-s)-2\psi(T-s)Q^{\top}Q\psi(t-s)\\
    &\quad + \frac{\gamma}{\gamma-1}(2\psi(T-s)Q^{\top}\rho v^{\top}+4\psi(T-s)Q^{\top}\rho\rho^{\top}Q\psi(T-s)).
    \end{aligned}
\end{equation*}
Let $L$ be the resolvent of the first kind of the integral kernel $K$.
Then by the associativity of the convolution and applying the fundamental theorem of calculus we obtain
\begin{equation*}
\begin{aligned}
f(\psi)(T-s) {}
    &= (\psi\ast L)'(T-s).
\end{aligned}
\end{equation*}
Thus we can write $F_s$ as
\begin{equation*}
    \begin{aligned}
    F_s{}
    & = (\psi\ast L)'(T-s)-\psi(T-s)M-M^{\top}\psi(T-s)-2\psi(T-s)Q^{\top}Q\psi(t-s)\\
    &\quad + \frac{\gamma}{\gamma-1}(2\psi(T-s)Q^{\top}\rho v^{\top}+4\psi(T-s)Q^{\top}\rho\rho^{\top}Q\psi(T-s)),
    \end{aligned}
\end{equation*}
We are interested in the expression $\int_0^T \operatorname{Tr}\big[F_s\Sigma_s\big] ds$.
Under the probability measure ${\Q}$, the process 
\[
\tilde{W}_{t}^{\sigma}= W_t^{\sigma}-\gamma\Sigma_t^{1/2}\pi_t^{*}\rho^{\top}
\]
is a $d\times d$-dimensional brownian motion by Girsanov's theorem. 
The dynamics of the volatility process $\Sigma$ under $\hat{\Q}$ can thus be written as
\begin{equation}
    \begin{aligned}
    \Sigma_t&
    =\Sigma_0+\int_0^t K(t-s)(NN^{\top}+M\Sigma_s+\Sigma_s M^{\top}+\frac{\gamma}{1-\gamma}\Sigma_s v\rho^{\top}Q\\
    &\quad +\frac{2\gamma}{1-\gamma}\Sigma_s\psi(T-s)Q^{\top}\rho\rho^{\top}Q+\frac{\gamma}{1-\gamma}Q^{\top}\rho v^{\top}\Sigma_s+\frac{2\gamma}{1-\gamma}Q^{\top}\rho\rho^{\top}Q\psi(T-s)\Sigma_s)ds\\
    &\quad + \int_0^t K(t-s)\Sigma_s^{1/2}d\tilde{W}_s^{\sigma}Q+Q^{\top}(d\tilde{W}_s^{\sigma})^{\top}\Sigma_s^{1/2}K(t-s).
    \end{aligned}
\end{equation}

We insert the dynamics of $\Sigma$ into the expression
$\operatorname{Tr}[\int_0^T (\psi\ast L)'(T-s)\Sigma_s ds]$ and simplify using calculus of convolutions and resolvents. Finally, some terms cancel out and we end up with
\begin{equation*}
    \begin{aligned}
    \int_0^T \operatorname{Tr}[F_s\Sigma_s]ds{}
    &= \operatorname{Tr}[\int_0^T f(\psi)(T-s)\Sigma_0 ds]+\operatorname{Tr}[\int_0^T \psi(T-s)NN^{\top} ds]\\ &\qquad-\operatorname{Tr}[\int_0^T 2\psi(T-s)Q^{\top}Q\psi(t-s)\Sigma_s ds]\\ 
    &\qquad +\operatorname{Tr}[\int_0^T \psi(T-s)( \Sigma_s^{1/2}d\tilde{W}_s^{\sigma}Q+Q^{\top}(d\tilde{W}_s^{\sigma})^{\top}\Sigma_s^{1/2})].
    \end{aligned}
\end{equation*}
Hence we obtain
\begin{equation*}
    \begin{aligned}
    {}
    & x_0^{-\gamma}\EX_{x_0,\Sigma_0}[{(X_T^{\pi^*})^{\gamma}}]=\EX_{x_0,\Sigma_0}^{\Q}\Big[\expo\Big(\int_0^T \gamma r_s ds+\int_0^T \operatorname{Tr}[F_s V_s] ds\Big)\Big]\\
    &= \expo\Big(\int_0^T \gamma r_s ds + \operatorname{Tr}[\int_0^T f(\psi)(T-s)\Sigma_0 ds]+\operatorname{Tr}[\int_0^T \psi(T-s)NN^{\top} ds]\Big)\\
    &\quad \times\EX_{x_0,\Sigma_0}^{\Q}\Big[\expo\Big(-2\operatorname{Tr}[\int_0^T Q\psi(T-s)\Sigma_s\psi(t-s)Q^{\top} ds]  + 2\operatorname{Tr}[\int_0^T Q\psi(T-s)\Sigma_s^{1/2}d\tilde{W}_s^{\sigma}] \Big)\Big].\\
    \end{aligned}
\end{equation*}
Since $\psi$ is continuous, it is bounded and therefore the stochastic exponential is a true $\Q$- martingale with expectation $1$ by Lemma \ref{A1}. Thus we get
\[
\EX_{x_0,\Sigma_0}\Big[{\frac{1}{\gamma}(X_T^{\pi^*})^{\gamma}}\Big]=\frac{x_0^\gamma}{\gamma}\expo\Big(\int_0^T \gamma r_s + \operatorname{Tr}[f(\psi)(T-s)\Sigma_0 +\psi(T-s)NN^{\top}]ds\Big).
\]
This completes the first part of the proof.\\
\textbf{Proof of the the inequality \eqref{b1}}
{First note that standard techniques to prove the inequality like \cite[Proposition 4.5]{Kra05} do not apply in our setting, since we can not use Itô's formula due to the non-Markovianity of our volatility model. Therefore, we use a different approach writing an arbitrary strategy $\pi$ in terms of $\pi^*$ and some remainder $\hat{\pi}$.}
To this end, let $\pi_t$ be an arbitrary admissible portfolio strategy. Define $\hat{\pi}:=\pi-\pi^*$ and write $\pi_t=\pi_t^* +\hat{\pi}_t$, where $\pi_t^*=\frac{1}{1-\gamma}(v+2\psi(T-t)Q^{\top}\rho)$. 
Since $\pi_t$ is bounded by assumption, we can define a new probability measure $\hat{\Q}$ with Radon-Nikodym density  
\begin{equation}\label{density hatQ}
\hat{Z}_0:=\frac{d\hat{\Q}}{d\mathbb{P}}|_{\mathcal{F}_0}=\expo\Big(\gamma\int_0^T \pi_s^{\top}\Sigma_s^{1/2}dW_s^{S}-\frac{\gamma^2}{2}\int_0^T \lVert \pi_s^{\top}\Sigma_s^{1/2} \rVert_{2}^2 ds\Big)
\end{equation}
by Lemma \ref{A2}.   
Analogously to the first part we obtain
\begin{equation*}
    \begin{aligned}
    & x_0^{-\gamma}\EX_{x_0,\Sigma_0}{(X_T^{\pi})^{\gamma}}=\EX_{x_0,\Sigma_0}^{\hat{\Q}}\Big[\expo\Big(\int_0^T \gamma r_s ds+\int_0^T \operatorname{Tr}\big[\hat{F}_s\Sigma_s\big] ds\Big)\Big].
    \end{aligned}
\end{equation*}
with
\begin{equation*}
\begin{aligned}
\hat{F}_s {}
    &:= \gamma v\pi_s^{\top}
    +\frac{\gamma(\gamma-1)}{2} \pi_s\pi_s^{\top} \\
    &= F_s + \gamma v \hat{\pi}_s^{\top} + \gamma(\gamma-1)\pi_s^*\hat{\pi}_s^{\top} + \frac{\gamma(\gamma-1)}{2}\hat{\pi}_s\hat{\pi}_s^{\top},
\end{aligned}
\end{equation*}
and $F_s$ is like in the first part. 
Under the probability measure $\hat{\Q}$ defined in \eqref{density hatQ}, the process 
\[
\hat{W}_{t}^{\sigma}=W_t^{\sigma}-\gamma\Sigma_t^{1/2}\pi_t\rho^{\top}= W_t^{\sigma}-\gamma\Sigma_t^{1/2}\pi_t^{*}\rho^{\top}-\gamma\Sigma_t^{1/2}\hat{\pi}_t\rho^{\top}
\]
is a $d\times d$-dimensional brownian motion by Girsanov's theorem. 
Hence the dynamics of the volatility process $\Sigma$ under $\hat{\Q}$ can be written as
\begin{equation*}
    \begin{aligned}
    \Sigma_t&
    =\Sigma_0+\int_0^t K(t-s)(NN^{\top}+M\Sigma_s+\Sigma_s M^{\top}+\frac{\gamma}{1-\gamma}\Sigma_s v\rho^{\top}Q\\
    &\quad +\frac{2\gamma}{1-\gamma}\Sigma_s\psi(T-s)Q^{\top}\rho\rho^{\top}Q+\frac{\gamma}{1-\gamma}Q^{\top}\rho v^{\top}\Sigma_s+\frac{2\gamma}{1-\gamma}Q^{\top}\rho\rho^{\top}Q\psi(T-s)\Sigma_s\\
    &\quad + \gamma\Sigma_s\hat{\pi}\rho^{\top}Q+\gamma Q^{\top}\rho\hat{\pi}^{\top}\Sigma_s) + \int_0^t K(t-s)\Sigma_s^{1/2}d\hat{W}_s^{\sigma}Q+Q^{\top}(d\hat{W}_s^{\sigma})^{\top}\Sigma_s^{1/2}K(t-s)
    \end{aligned}
\end{equation*}
Again we calculate 
$\operatorname{Tr}[\int_0^T (\psi\ast L)'(T-s)\Sigma_s ds]$ inserting the dynamics of $\Sigma$ and using calculus of convolutions and resolvents. Some terms cancel out and this time we end up with
\begin{equation*}
\begin{aligned}
&\int_0^T \operatorname{Tr}[\hat{F}_s \Sigma_s] ds {}\\
    &\quad= \operatorname{Tr}[\int_0^T f(\psi)(T-s)\Sigma_0 ds]+\operatorname{Tr}[\int_0^T \psi(T-s)NN^{\top} ds]\\
    &\qquad-\operatorname{Tr}[\int_0^T 2\psi(T-s)Q^{\top}Q\psi(t-s)\Sigma_s ds] \\
    &\qquad +\operatorname{Tr}[\int_0^T \psi(T-s)( \Sigma_s^{1/2}d\hat{W}_s^{\sigma}Q+Q^{\top}(d\hat{W}_s^{\sigma})^{\top}\Sigma_s^{1/2})]+ \operatorname{Tr}[\int_0^T\frac{\gamma(\gamma-1)}{2}\hat{\pi}_s\hat{\pi}_s^{\top}\Sigma_s ds]
\end{aligned}
\end{equation*}
{The term $\operatorname{Tr}[\int_0^T\frac{\gamma(\gamma-1)}{2}\hat{\pi}_s\hat{\pi}_s^{\top}\Sigma_s ds]$ is equivalent to $\int_0^T\frac{\gamma(\gamma-1)}{2}\hat{\pi}_s\Sigma_s \hat{\pi}_s^{\top}ds$ and since $\Sigma_s$ is positive definite and $\gamma\in (0,1)$, it has to be less than or equal to $0$. }
Thus, for the expectation we get 
\begin{equation*}
    \begin{aligned}
    {}
    & x_0^{-\gamma}\EX_{x_0,\Sigma_0}[{(X_T^{\pi})^{\gamma}}]=\EX_{x_0,\Sigma_0}^{\hat{\Q}}\Big[\expo\Big(\int_0^T \gamma r_s ds+\int_0^T \operatorname{Tr}[\hat{F}_s\Sigma_s] ds\Big)\Big]\\
    &=\expo\Big(\int_0^T \gamma r_s ds + \operatorname{Tr}[\int_0^T f(\psi)(T-s)\Sigma_0 ds]+\operatorname{Tr}[\int_0^T \psi(T-s)NN^{\top} ds] \Big)\\
    &\qquad\times\EX_{x_0,\Sigma_0}^{\hat{\Q}}\Big[\expo\Big(-\operatorname{Tr}[\int_0^T 2\psi(T-s)Q^{\top}Q\psi(t-s)\Sigma_s ds]  + 2\operatorname{Tr}[\int_0^T Q\psi(T-s)\Sigma_s^{1/2}d\hat{W}_s^{\sigma}] \\
    & \qquad\qquad\qquad + \operatorname{Tr}[\int_0^T\frac{\gamma(\gamma-1)}{2}\hat{\pi}_s\Sigma_s\hat{\pi}_s^{\top} ds]\Big)\Big]\\
    &\leq \expo\Big(\int_0^T \gamma r_s ds + \operatorname{Tr}[\int_0^T f(\psi)(T-s)\Sigma_0 ds]+\operatorname{Tr}[\int_0^T \psi(T-s)NN^{\top} ds] \Big)\\ 
    &\quad \times\EX_{x_0,\Sigma_0}^{\hat{\Q}}\Big[\expo\Big(-2\operatorname{Tr}[\int_0^T Q\psi(T-s)\Sigma_s\psi(t-s)Q^{\top} ds]  + 2\operatorname{Tr}[\int_0^T Q\psi(T-s)\Sigma_s^{1/2}d\hat{W}_s^{\sigma}] \Big)\Big].
    \end{aligned}
\end{equation*}
{Note that the last inequality does not hold for the case $\gamma<0$}. Since the stochastic exponential is a $\hat{\Q}$-martingale with expectation $1$ by Lemma \ref{A1}, we finally obtain
\[
\EX_{x_0,\Sigma_0}\Big[{\frac{1}{\gamma}(X_T^{\pi})^{\gamma}}\Big]\leq \frac{x_0^\gamma}{\gamma}\expo\Big(\int_0^T \gamma r_s + \operatorname{Tr}[f(\psi)(T-s)\Sigma_0 +\psi(T-s)NN^{\top}]ds \Big),
\]
which completes the proof.} \qed

\section{Numerical examples}\label{sec: example}

In this section we compute the optimal portfolio strategy numerically in two-dimensional examples. To begin, we consider a financial market with one riskfree asset and $d=2$ risky assets and an investment horizon of $T=1$ year. The parameters are taken from \cite{BPT10}, where such a model is calibrated to real market data from the Standard and Poor’s 500 Index and 30-year Treasury bond. 
They obtained the following estimation for the model parameters:

\[   
    M=\left( \begin{array}{cc}
-1.21 & 0.491 \\
0.3292 & -1.271
\end{array} \right),\qquad
Q=\left( \begin{array}{cc}
0.167 & 0.033 \\
0.001 & 0.09
\end{array} \right),
\]

  \[ \rho=\left( \begin{array}{cc}
-0.115 \\
-0.549 
\end{array} \right),\qquad
v=\left( \begin{array}{cc}
4.722  \\
3.317
\end{array} \right),
\]
and $NN^{\top}={10} Q^{\top}Q$.
Roughness of the model is obtained by taking an appropriate integration kernel.
We choose a fractional kernel of the form
\[   
    K(t)=\left( \begin{array}{cc}
\frac{t^{\alpha-1}}{\Gamma(\alpha)} & 0 \\
0 & \frac{t^{\alpha-1}}{\Gamma(\alpha)}
\end{array} \right),\qquad \alpha\in(\frac{1}{2},1).
\]
This corresponds to the rough Heston model in the one-dimensional case. The roughness of the volatility paths is determined by the parameter $\alpha$ and for $\alpha\rightarrow 1$ we recover the classical model. The parameter $\alpha$ is linked to the Hurst parameter $H$ via the equation $\alpha=H+\frac{1}{2}$. In our example the investor has a power utility function
    \[
    U(x)=\frac{1}{\gamma}x^{\gamma},\quad 0<\gamma<1. 
    \]
The optimal strategy $\pi^*$ consists of a constant term $$\frac{v}{1-\gamma}$$ and a time-depending term $$\frac{2}{1-\gamma}\psi(T-t)Q^{\top}\rho,$$ the hedging demand. Here the parameter $v$ is the market price of risk, $\gamma$ is the relative risk aversion, $\rho$ is the correlation vector, $Q$ is the Matrix governing the diffusion of the volatility process and $\psi$ is the solution of the matrix-valued Riccati-Volterra equation \eqref{MRV}. 
In order to compute the hedging demand, we have to find the solution $\psi$ of equation \eqref{MRV}. To this end we use the fractional Adams method developed in \cite{DFF02},\cite{DFF04} to obtain a numerical solution.     
The next diagram shows that if the roughness level $\alpha\rightarrow 1$, we recover the results of \cite[Figure 1]{BL13} for the classical Wishart model.     
    
\begin{figure}[!ht]
\centering
    \begin{subfigure}{0.45\linewidth}
        \includegraphics[width=\linewidth]{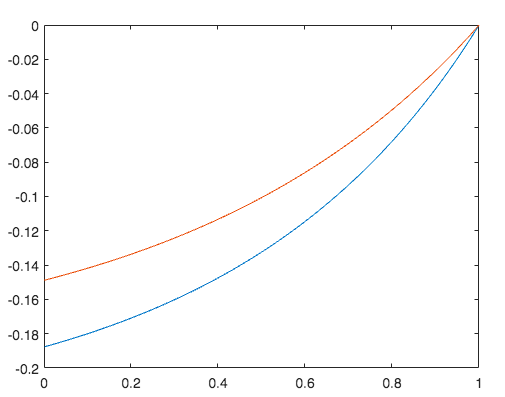}
    \caption{}
    \end{subfigure}
\hfil
    \begin{subfigure}{0.45\linewidth}
        \includegraphics[width=\linewidth]{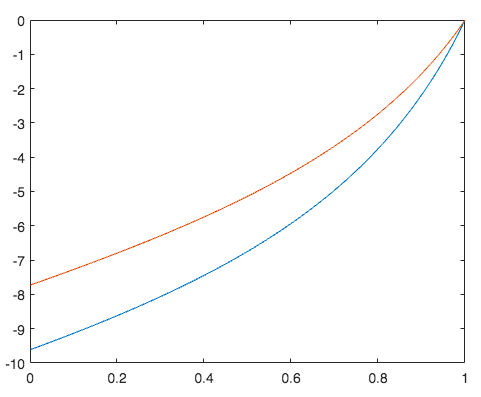}
    \caption{}
    \end{subfigure}

\caption{Hedging demands for roughness level $\alpha=1$ for parameter $\gamma=0.2$ (A) and $\gamma=0.8$ (B).}
    \label{fig:my figure2}
    \end{figure}        

\begin{figure}[!ht]
\centering
    \begin{subfigure}{0.45\linewidth}
        \includegraphics[width=\linewidth]{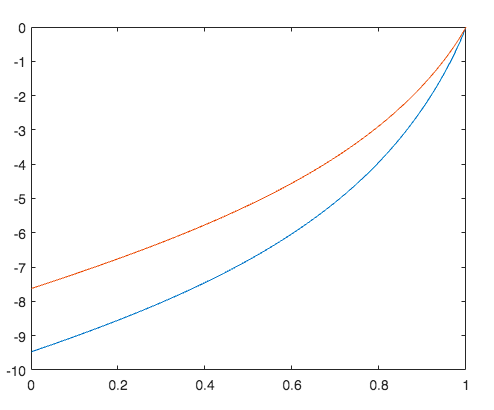}
    \caption{}
    \end{subfigure}
\hfil
    \begin{subfigure}{0.45\linewidth}
        \includegraphics[width=\linewidth]{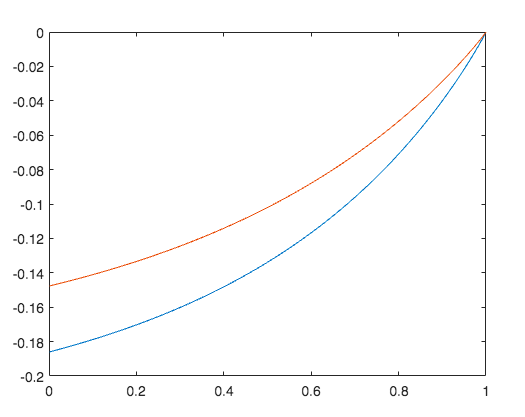}
    \caption{}
    \end{subfigure}

    \begin{subfigure}{0.45\linewidth}
        \includegraphics[width=\linewidth]{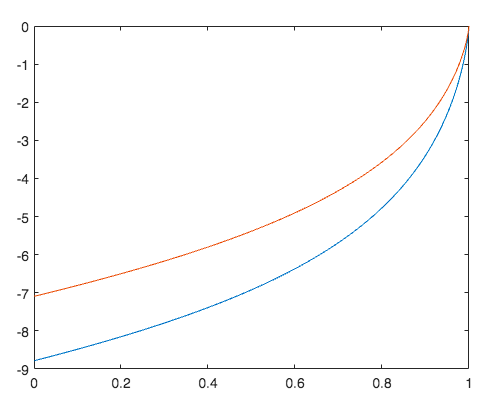}
    \caption{}
    \end{subfigure}
\hfil
    \begin{subfigure}{0.45\linewidth}
        \includegraphics[width=\linewidth]{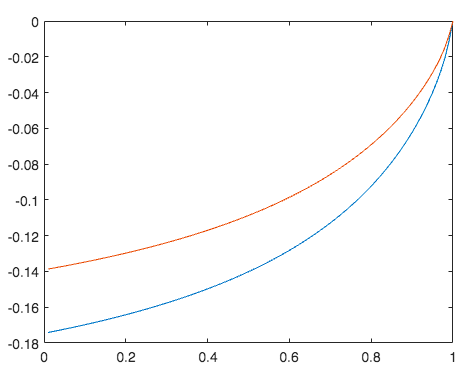}
    \caption{}
    \end{subfigure}

    \begin{subfigure}{0.45\linewidth}
        \includegraphics[width=\linewidth]{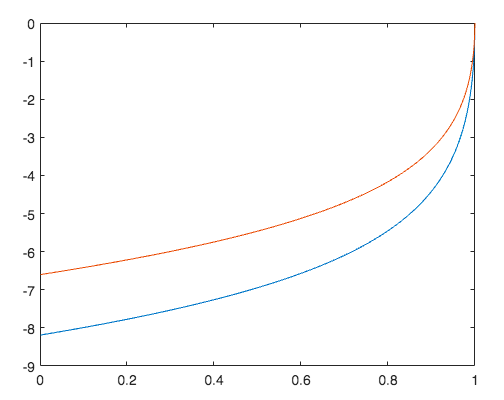}
    \caption{}
    \end{subfigure}
\hfil
    \begin{subfigure}{0.45\linewidth}
        \includegraphics[width=\linewidth]{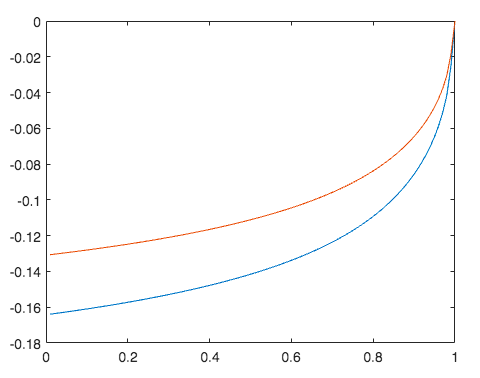}
    \caption{}
    \end{subfigure}
\caption{Left hand side: Hedging demands $\frac{2}{1-\gamma}\psi(T-t)Q^{\top}\rho$ for risk aversion parameter $\gamma=0.8$ and roughness levels $\alpha=0.95$ (A), $\alpha=0.75$ (C) and $\alpha=0.55$ (E).\newline\newline
Right hand side: Hedging demands $\frac{2}{1-\gamma}\psi(T-t)Q^{\top}\rho$ for risk aversion parameter $\gamma=0.2$ and roughness levels $\alpha=0.95$ (B), $\alpha=0.75$ (D) and $\alpha=0.55$ (F).}
    \label{fig:my figure}
    \end{figure}
    

%

Figure \ref{fig:my figure} shows that, in accordance with \cite{BL13}, the hedging demand for $\gamma\in (0,1)$ is negative. The lower the risk aversion of the investor, the more negative is his hedging demand. The roughness of the volatility of the assets also affects the hedging demand over time. Our illustrations show that the curvature of the hedging demand is increasing as the paths of the volatility become rougher. {This behaviour is in line with the results of \cite{HW20} in the case of one risky asset.}          

\begin{figure}[!ht]
\centering
    \begin{subfigure}{0.45\linewidth}
        \includegraphics[width=\linewidth]{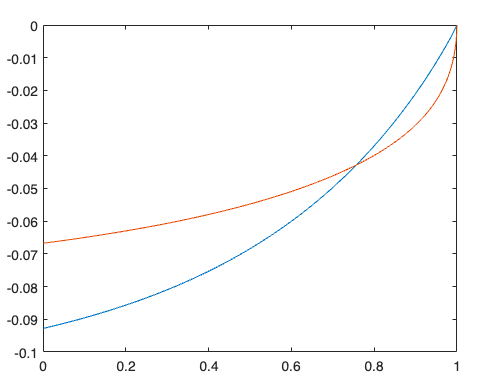}
    \caption{Vol of vol $\frac{1}{2}Q$}
    \end{subfigure}
\hfil
    \begin{subfigure}{0.45\linewidth}
        \includegraphics[width=\linewidth]{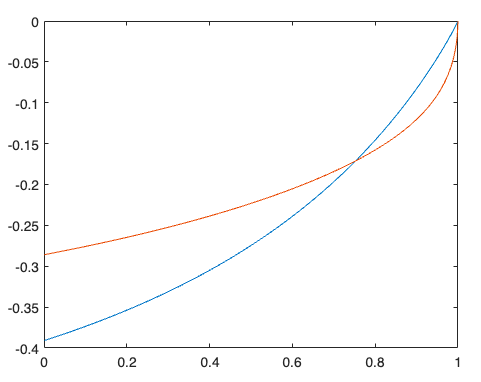}
   \caption{Vol of vol $2Q$}
   \end{subfigure}
\caption{Effect of two different levels of the volatility of volatility on the hedging demand. The hedging demand depends linearly on the vol of vol $Q$.}
\label{fig:my figure5}
\end{figure}

\begin{figure}[!ht]
\centering
    \begin{subfigure}{0.45\linewidth}
        \includegraphics[width=\linewidth]{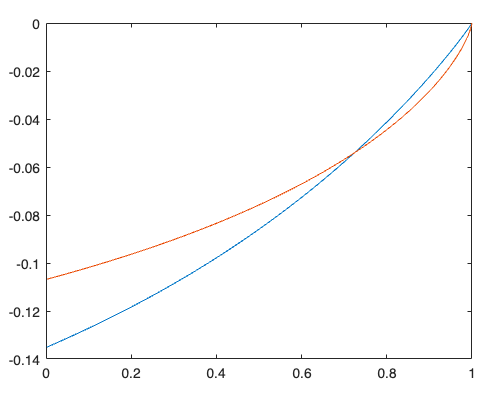}
    \caption{$T=0.5,\ \alpha_1=0.95,\ \alpha_2=0.75$}
    \end{subfigure}
\hfil
    \begin{subfigure}{0.45\linewidth}
        \includegraphics[width=\linewidth]{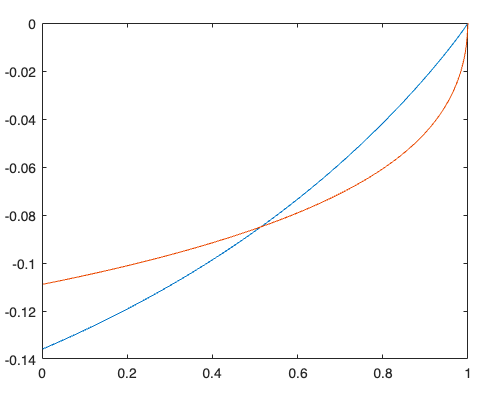}
    \caption{$T=0.5,\ \alpha_1=0.95,\ \alpha_2=0.55$}
    \end{subfigure}

    \begin{subfigure}{0.45\linewidth}
        \includegraphics[width=\linewidth]{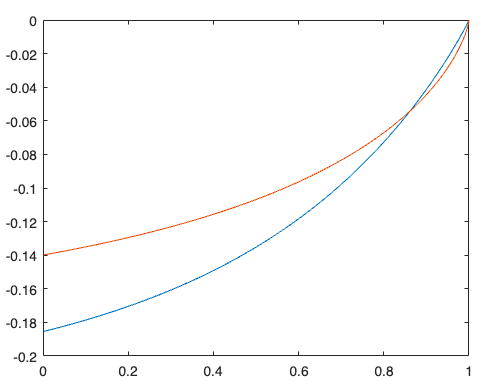}
    \caption{$T=1,\ \alpha_1=0.95,\ \alpha_2=0.75$}
    \end{subfigure}
\hfil
    \begin{subfigure}{0.45\linewidth}
        \includegraphics[width=\linewidth]{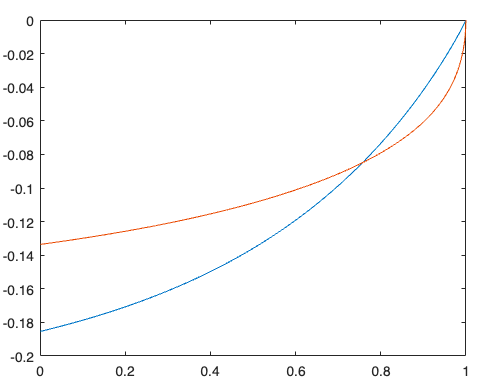}
    \caption{$T=1,\ \alpha_1=0.95,\ \alpha_2=0.55$}
    \end{subfigure}

    \begin{subfigure}{0.45\linewidth}
        \includegraphics[width=\linewidth]{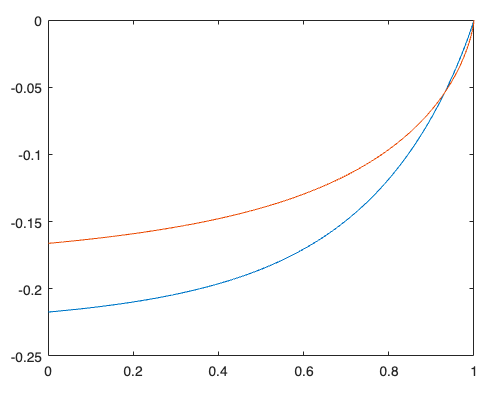}
    \caption{$T=2,\ \alpha_1=0.95,\ \alpha_2=0.75$}
    \end{subfigure}
\hfil
    \begin{subfigure}{0.45\linewidth}
        \includegraphics[width=\linewidth]{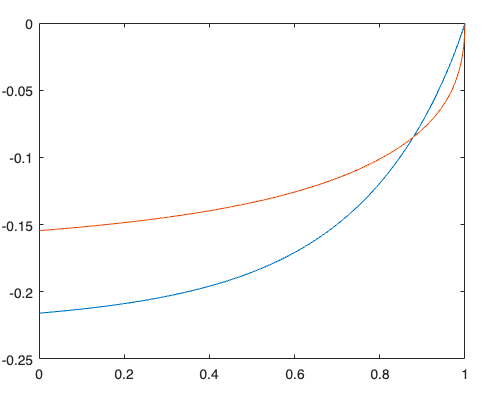}
    \caption{$T=2,\ \alpha_1=0.95,\ \alpha_2=0.55$}
    \end{subfigure}
\caption{Effect of the time horizon $T$ on the hedging demand for different levels of roughness for an investor with risk aversion $\gamma=0.2$.}
\label{fig:my figure3}
\end{figure}

\begin{figure}[!ht]
\centering
    \begin{subfigure}{0.45\linewidth}
        \includegraphics[width=\linewidth]{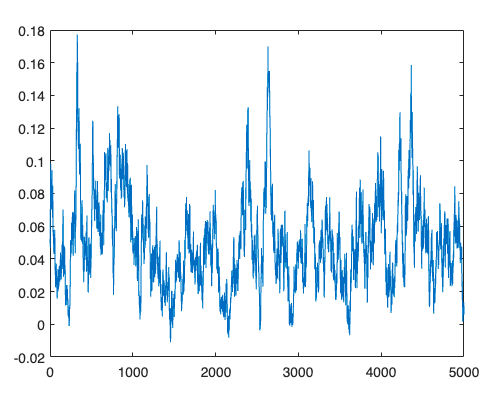}
    \caption{sample path for positive asset correlation}
    \end{subfigure}
\hfil
    \begin{subfigure}{0.45\linewidth}
        \includegraphics[width=\linewidth]{0.95,0.55,1.png}
    \caption{$T=1,\ \alpha_1=0.95,\ \alpha_2=0.55$}
    \end{subfigure}

    \begin{subfigure}{0.45\linewidth}
        \includegraphics[width=\linewidth]{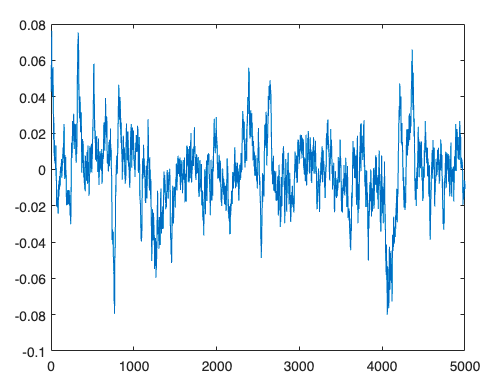}
    \caption{sample path for asset correlation around $0$}
    \end{subfigure}
\hfil
    \begin{subfigure}{0.45\linewidth}
        \includegraphics[width=\linewidth]{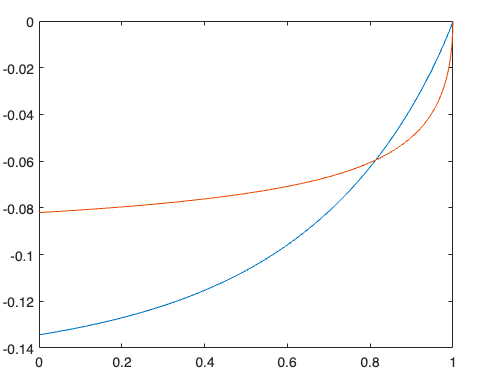}
    \caption{$T=1,\ \alpha_1=0.95,\ \alpha_2=0.55$}
    \end{subfigure}
    
  \begin{subfigure}{0.45\linewidth}
        \includegraphics[width=\linewidth]{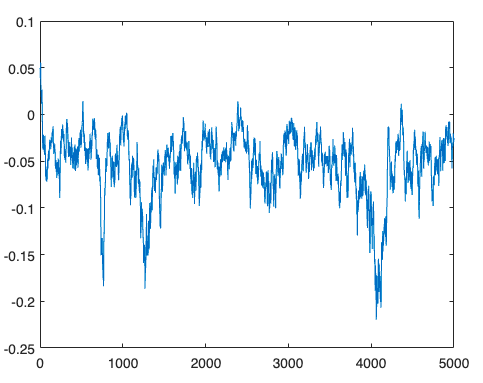}
    \caption{sample path for negative asset correlation}
    \end{subfigure}
\hfil
    \begin{subfigure}{0.45\linewidth}
        \includegraphics[width=\linewidth]{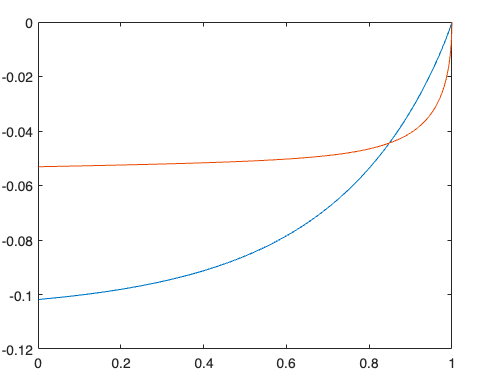}
    \caption{$T=1,\ \alpha_1=0.95,\ \alpha_2=0.55$}
    \end{subfigure}

\caption{Effect of the asset correlation on the hedging demand.}
\label{fig:my figure4}
\end{figure}

{
To have a comparison with another multivariate Volterra model, we adapt a numerical experiment from \cite[Chapter 6]{PJM21} which was carried out for the Markowitz problem in a rough Stein-Stein model. Our aim is to investigate how the optimal investment strategy is affected on different time intervals within the investment horizon, if the investor can choose among a rough and a smooth asset with Hurst parameters $H_1<H_2$. For this experiment we again use the parameters $M,Q,\rho$ and $v$ as in the beginning of this section. Recall that the Hurst parameter $H$ is linked to the parameter $\alpha$ in our integral kernels via {$H=\alpha-\frac{1}{2}$}. 
It turns out that as in \cite{PJM21}, we end up with three distinct regimes:
\begin{enumerate}
    \item $T<<1$:  When the investment horizon is close to the end, the investor is selling a larger amount of the smooth asset than of the rough one.
    \item $T\sim 1$: Here a transition of the investors behavior appears. First the agent prefers selling the rough asset but as the final horizon approaches, he prefers to sell the smooth asset.
    \item $T>>1$: All the time until the transition point close to the maturity, the investor prefers selling the rough asset.
\end{enumerate}
Our results are illustrated in Figure \ref{fig:my figure3}. As pointed out in \cite{PJM21}, a possible interpretation of this transition could be that rough processes are more volatile than smooth processes in the short term but less volatile in the long term. Thus, when there is not much time left, the investor prefers rough assets to obtain some performance, whereas he favors the smooth asset on the long run. It turns out that the larger the difference of the roughness of the two assets becomes, the more the effect described above is pronounced. This means that if the second asset is very rough, the transition of the investors behavior happens earlier and the amount of rough assets the investor is selling decreases faster.  
In contrast to the Markowitz portfolio strategies studied in \cite{PJM21}, we do not observe a structural difference in our portfolio strategies for different levels of asset correlation. While a positive asset correlation leads to a \textit{buy rough sell smooth} strategy (cf. \cite{GH20}) for the optimization problem in \cite{PJM21}, in our case the correlation level only affects the extent to which the investor is selling both the smooth and the rough asset (see Figure \ref{fig:my figure4}). We want to point out that an average asset correlation of $0$ can be obtained by setting the non-diagonal entries of the matrices $M$ and $Q$ to zero. 
Since our hedging demands are linear in terms of the vol-of-vol $Q$, the investors' preferences are preserved in our model if $Q$ is multiplied with a constant, whereas in [AJMP21] the investor's preferences do depend on the vol-of-vol. For completeness, this is illustrated in Figure~\ref{fig:my figure5}.}

\newpage

\newpage    

\newpage



\section{Existence and uniqueness results for the Volterra equations}\label{Existence}

In this section we provide existence and uniqueness results for the Riccati-Volterra equations \eqref{MRV} and the stochastic Volterra equation \eqref{vola}. 

\subsection{Solution of the Volterra Wishart process}
Previous research on stochastic Volterra equations has been carried out in the vector framework (cf. \cite{PJM21},\cite{AJ19}).
For proving existence of a weak solution for the matrix-valued Volterra-Wishart equation \eqref{vola}, we use the vectorization operator in order to be able to resort to existing literature.

\begin{thm}\label{ExVol}
Assume that $K$ admits a resolvent of the first kind and that the components of $K$ satisfy \eqref{kernels}. Then the stochastic Volterra equation \eqref{vola} has a unique in law $\R^{d\times d}$-valued continuous weak solution on the interval $[0,T_{\operatorname{pos}})$,
where $T_{\operatorname{pos}}:={\operatorname{sup}}\{t\geq 0: \tau\geq t \text{ almost sure}\}$ $\tau:=\operatorname{inf}\{t\geq 0: \Sigma_t \text{ is not positive definite}\}$.
\end{thm}

\emph{Proof:}
We cast the problem into the vector framework using the vectorization operator $\operatorname{vec}: \R^{d\times d}\rightarrow \R^{d\cdot d}$ stacking the columns of a matrix on top of one another. This leads to the $\R^{d\cdot d}$-valued stochastic Volterra equation

\begin{equation*}
    \begin{aligned}
    \operatorname{vec}(\Sigma){}
    & = \vecop(K\ast (NN^{\top}+M\Sigma+\Sigma M^{\top}+\Sigma^{1/2}dW^{\sigma}Q+Q^{\top}(dW^{\sigma})^{\top}\Sigma^{1/2}))\\
    & =(I\otimes K)\ast\vecop(NN^{\top}+M\Sigma+\Sigma M^{\top}+\Sigma^{1/2}dW^{\sigma}Q+Q^{\top}(dW^{\sigma})^{\top}\Sigma^{1/2})\\
    &= (I\otimes K)\ast[\vecop(NN^{\top})+ (I\otimes M + M\otimes I)\vecop(\Sigma)\\
    &\qquad\qquad\qquad\qquad\qquad + Q^{\top}\otimes \Sigma^{1/2}\vecop(dW^{\sigma})+\Sigma^{1/2}\otimes Q^{\top}\vecop((dW^{\sigma})^{\top})]\\
    & = (I\otimes K)\ast[b(\vecop(\Sigma))dt + c(\vecop(\Sigma))d\vecop(W^{\sigma})]
    \end{aligned}
\end{equation*}
where 
\[
b(\vecop(\Sigma)):=\vecop(NN^{\top})+ (I\otimes M + M\otimes I)\vecop(\Sigma)
\]
and 
\[
c(\vecop(\Sigma)):=Q^{\top}\otimes \Sigma^{1/2}+ f(\Sigma^{1/2}\otimes Q^{\top}),
\]
{and $f: \R^{d\times d} \rightarrow \R^{d\times d}$ swaps entries such that $f(\Sigma^{1/2}\otimes Q^{\top})d\operatorname{vec}(W^{\sigma})=\Sigma^{1/2}\otimes Q^{\top}d\operatorname{vec}((W^{\sigma})^{\top})$}. We now show that $b$ and $c$ fulfill a linear growth condition.

\begin{equation*}
    \begin{aligned}
        \abs{b(\vecop(\Sigma))}
        & \leq \abs{\vecop(NN^{\top})}+\lVert I\otimes M + M\otimes I \rVert_F \abs{\vecop(\Sigma)}\\
        & \leq  \abs{\vecop(NN^{\top})}+ 2\sqrt{d}\lVert M \rVert_F \abs{\vecop(\Sigma)}\leq C_1 (1+\abs{\vecop(\Sigma)})
    \end{aligned}
\end{equation*}
with $C_1=\operatorname{max}(\abs{\vecop(NN^{\top})}, 2\sqrt{d}\lVert M \rVert_F)$.

\begin{equation*}
    \begin{aligned}
        \abs{c(\vecop(\Sigma))}
        & \leq \lVert Q^{\top}\otimes \Sigma^{1/2}\rVert_F + \lVert f(\Sigma^{1/2}\otimes Q^{\top})\rVert_F = 2\lVert Q^{\top} \rVert_F \lVert \Sigma^{1/2}\rVert_F
    \end{aligned}
\end{equation*}

Using the fact that $\lVert \Sigma^{1/2}\rVert_F = \sqrt{\operatorname{Tr}(\Sigma^{1/2} \Sigma^{1/2})}$, we obtain
\begin{equation*}
    \begin{aligned}
        \abs{c(\vecop(\Sigma))}^2
        &  \leq 4\lVert Q^{\top} \rVert_F^2 \lVert \Sigma^{1/2}\rVert_F^2 =  4\lVert Q^{\top} \rVert_F^2 \operatorname{Tr}(\Sigma^{1/2} \Sigma^{1/2}) = 4\lVert Q^{\top} \rVert_F^2 \operatorname{Tr}(\Sigma)\\
        & \leq 4\lVert Q^{\top} \rVert_F^2 \sqrt{d}(1+\lVert \Sigma \rVert_F) \leq  4\lVert Q^{\top} \rVert_F^2 \sqrt{d} (1+\lVert \Sigma \rVert_F)^2 = C_2^2 (1+\abs{\vecop(\Sigma)})^2.
    \end{aligned}
\end{equation*}
Here we have used the fact that 
\begin{equation*}
    \begin{aligned}
        \operatorname{Tr}(\Sigma)
        &= \sum\limits_{i=1}^d \sigma_{ii}\leq \sqrt{d}\sqrt{\sum\limits_{i=1}^d \sigma_{ii}^2}\leq \sqrt{d}(1+\sum\limits_{i=1}^d \sigma_{ii}^2)\leq \sqrt{d}(1+\sum\limits_{i,j=1}^d \sigma_{ij}^2)=\sqrt{d}(1+\operatorname{Tr}(\Sigma^2)).
    \end{aligned}
\end{equation*}
Hence for $C_2=2 d^{1/4} \lVert Q \rVert_F$, we get
\[
\abs{c(\vecop(\Sigma))}\leq C_2 (1+\abs{\vecop(\Sigma)}).
\]
Thus $b$ and $c$ fulfill the linear growth condition \cite[Condition 3.1]{AJ19} and hence we can apply \cite[Theorem 3.4]{AJ19} to get the desired result. \qed\\

\subsection{Solution of the Riccati Volterra equation}

We use the concept of non-continuable solutions described in \cite{AJ19}.
For an integral kernel $K_i\in L_{\operatorname{loc}}^2(\R_+,\R^{d\times d})$ and a function $f: \R^d\rightarrow \R^d$ consider the Volterra integral equation
\begin{equation}\label{Volterra equ}
    \psi = K\ast f(\psi).
\end{equation}
For $0<T_{\operatorname{max}}<\infty$, by a non-continuable solution of \eqref{Volterra equ} we denote a pair $(\psi,T_{\operatorname{max}})$ with $\lVert \psi \rVert_{L_2[0,T_{\operatorname{max}})}<\infty$, such that $\psi$ satisfies \eqref{Volterra equ} on $[0,T_{\operatorname{max}})$ and $\lVert \psi \rVert_{L_2[0,T_{\operatorname{max}}]}=\infty$. A non-continuable solution $(\psi,T_{\operatorname{max}})$ is unique if for any $T>0$ and $\tilde{\psi}$ with $\lVert \tilde{\psi} \rVert_{L_2[0,T]}<\infty$ satisfying \eqref{Volterra equ} on $[0,T]$, we have $T<T_{\operatorname{max}}$ and $\tilde{\psi}=\psi$ on $[0,T]$. If $T_{\operatorname{max}}=\infty$, we call $\psi$ a global solution of \eqref{Volterra equ}.

\begin{thm}\label{vecRicEx}
For each of the Riccati-Volterra equations \eqref{RV1} and \eqref{RV2} there exists a unique non-continuable continuous solution $(\psi, T_{\operatorname{max}})$. 
\end{thm}

\emph{Proof:}
The Riccati-Volterra equations \eqref{RV1} and \eqref{RV2} are both of the form 
\[
\chi = F(\chi)\ast K
\]
where $F$ is given by
\[
F(\chi)=a+ B\chi +(c_1\chi_1^2,\dots, c_d \chi_d^2)^{\top}
\]
with $a\in \R^d$, $B\in \R^{d\times d}$ and $c_i\in \R$.
Thus for $x$, $y\in \R^d$ we get

\begin{equation*}
    \begin{aligned}
   \abs{F(x)-F(y)}{}
    &= \abs{Bx +(c_1 x_1^2,\dots, c_d x_d^2)^{\top}-By - (c_1 y_1^2,\dots, c_d y_d^2)^{\top} }\\
    &\leq \abs{B(x-y)}+\abs{(c_1 (x_1^2-y_1^2),\dots, c_d (x_d^2-y_d^2))^{\top}}\\
    &\leq \abs{x-y}\lVert B \rVert_F + \max_{1\leq i\leq d}(\abs{c_i})\abs{((x_1+y_1)(x_1-y_1),\dots, (x_d+y_d)(x_d-y_d))^{\top}} \\
    &\leq  \abs{x-y}\lVert B \rVert_F + \max_{1\leq i\leq d}(\abs{c_i})\abs{x+y}\cdot \abs{x-y} \\
    &\leq C_1 \abs{x-y} + C_2 \abs{x-y}(\abs{x}+\abs{y})
    \end{aligned}
\end{equation*}
with positive constants $C_1=\lVert B \rVert_F$ and $C_2=\max_{1\leq i\leq d}(\abs{c_i})$. Here we have used the triangle inequality, the Cauchy-Schwarz inequality and the fact that for $\alpha$, $\beta\in {\R_+^d}$ we get $\sum_{i=1}^n \alpha_i\beta_i\leq \sum_{i=1}^n \alpha_i \sum_{j=1}^n \beta_j$. Now existence and uniqueness of a non-continuable solution $(\chi, T_{\operatorname{max}})$ follows directly from \cite[Theorem B1]{AJ19}. \qed\\

In case that the matrix $D$ in equation \eqref{Vola} is diagonal, there exists even a global solution. 

\begin{thm}
Assume that the matrix $D$ in the drift of the volatility process is a diagonal matrix, i.e. $D=\diag{(\delta_1,\dots, \delta_d)}$. Then the Riccati-Volterra equation \eqref{RV2} has a unique global solution if for all $1\leq i\leq d$ 
\[
\delta_i+\frac{\gamma}{1-\gamma}\nu_i\rho_i\theta_i <0\text{ and }  (\delta_i+\frac{\gamma}{1-\gamma}\nu_i\rho_i\theta_i)^2-\frac{\gamma}{1-\gamma}(\frac{1-\gamma+\gamma\rho_i^2}{1-\gamma })\nu_i^2\theta_i^2 >0.
\]
\end{thm}

\emph{Proof:}
If $D$ is a diagonal matrix, the matrix $\Lambda$ in equation \eqref{RV2} becomes a diagonal matrix of the form $\Lambda=-\diag(\lambda_1,\dots,\lambda_d)$ with $\lambda_i=-\delta_i-\frac{\gamma}{1-\gamma}\nu_i\rho_i\theta_i$. Now the vector valued equation \eqref{RV2} can be decomposed into $d$ real valued Riccati-Volterra equations such that for the $i$th component $\psi_i$ of $\psi$ we obtain 
\begin{equation}\label{comp}
\psi_i(t)=\int_0^t K_i(t-s)\Big[ \frac{\gamma}{2(1-\gamma)}\theta_i^2-\lambda_i\psi_i(s)+\frac{1}{2}\nu_i^2\frac{1-\gamma+\gamma\rho_i^2}{1-\gamma}\psi_i^2(s)\Big]ds
\end{equation}
By our assumptions $\lambda_i>0$ and $\lambda_i^2-2\frac{\gamma \theta_i^2}{2(1-\gamma)}\nu_i^2\frac{1-\gamma+\gamma\rho_i^2}{1-\gamma}>0$ and therefore \cite[Lemma A2]{HW20} (cf. \cite[Lemma A5]{GKR19}) guarantees the existence of a unique continuous global solution of the equation \eqref{comp}. Combining the component-wise solutions, we finally obtain the unique global solution $\psi$ of equation \eqref{RV2}. \qed

A vectorization argument allows us to prove existence of a local solution for the matrix Riccati-Volterra equation \eqref{MRV}.

\begin{thm}\label{ExRic}
The matrix Riccati-Volterra equation \eqref{MRV} has a unique non-continuable continuous solution $(\psi,T_{\text{max}})$.
\end{thm}

\emph{Proof:}
We cast the problem into the vector framework using the vectorization operator $\operatorname{vec}: \R^{d\times d}\rightarrow \R^{d\cdot d}$ stacking the columns of a matrix on top of one another. This leads to the $\R^{d\cdot d}$-valued Riccati-Volterra equation

\begin{equation*}
    \begin{aligned}
    \operatorname{vec}(\psi){}
    & = \operatorname{vec}(K\ast[\psi\tilde{M}+\tilde{M}^{\top}\psi+2\psi\tilde{Q}^{\top}\tilde{Q}\psi+\tilde{\Gamma}])\\
    & = (I\otimes K)\ast \operatorname{vec}(\psi\tilde{M}+\tilde{M}^{\top}\psi+2\psi\tilde{Q}^{\top}\tilde{Q}\psi+\tilde{\Gamma}) =: (I\otimes K)\ast p(\operatorname{vec}(\psi)),
    \end{aligned}
\end{equation*}
where $\otimes: \R^{d\times d}\times \R^{d\times d}\rightarrow \R^{d\cdot d\times d\cdot d}$ denotes the Kronecker product. We now show that $p$ fulfills the growth condition 
\[
\abs{p(\operatorname{vec}(X))-p(\operatorname{vec}(Y))}\leq C_1 \abs{\operatorname{vec}(X)-\operatorname{vec}(Y)}+C_2 \abs{\operatorname{vec}(X)-\operatorname{vec}(Y)}(\abs{\operatorname{vec}(X)}+\abs{\operatorname{vec}(Y)})
\]
with positive constants $C_1$ and $C_2$.

\begin{equation*}
    \begin{aligned}
     &\abs{p(\operatorname{vec}(X))-p(\operatorname{vec}(Y))} {}\\
    & = \abs{\operatorname{vec}(X\tilde{M}+\tilde{M}^{\top}X+2X\tilde{Q}^{\top}\tilde{Q}X+\tilde{\Gamma})- \operatorname{vec}(Y\tilde{M}+\tilde{M}^{\top}Y+2Y\tilde{Q}^{\top}\tilde{Q}Y+\tilde{\Gamma})}\\
    & = \abs{\vecop((X-Y)\tilde{M})+\vecop(\tilde{M}^{\top}(X-Y))+\vecop(2X\tilde{Q}^{\top}\tilde{Q}X-2Y\tilde{Q}^{\top}\tilde{Q}Y)}\\
    &\leq \abs{\vecop((X-Y)\tilde{M})}+\abs{\vecop(\tilde{M}^{\top}(X-Y))}+2\abs{\vecop(X\tilde{Q}^{\top}\tilde{Q}X-Y\tilde{Q}^{\top}\tilde{Q}Y)}
    \end{aligned}
\end{equation*}

For the first and the second term of this sum, we have

\[
\abs{\vecop((X-Y)\tilde{M})}= \abs{(\tilde{M}^{\top}\otimes I)\vecop(X-Y)}\leq \lVert \tilde{M}^{\top}\otimes I\rVert_F \abs{\vecop(X-Y)}
\]
and
\[
\abs{\vecop(\tilde{M}^{\top}(X-Y))}= \abs{(I\otimes\tilde{M}^{\top})\vecop(X-Y)}\leq \lVert I\otimes\tilde{M}^{\top}\rVert_F \abs{\vecop(X-Y)}.
\]
Since $\lVert \tilde{M}^{\top}\otimes I\rVert_F= \lVert I\otimes\tilde{M}^{\top}\rVert_F =\sqrt{d}\lVert \tilde{M}\rVert_F $, we choose $C_1=2\sqrt{d}\lVert \tilde{M}\rVert_F$.
For the third term, it holds that
\[
\abs{\vecop(X\tilde{Q}^{\top}\tilde{Q}X-Y\tilde{Q}^{\top}\tilde{Q}Y)}=\lVert X\tilde{Q}^{\top}\tilde{Q}X-Y\tilde{Q}^{\top}\tilde{Q}Y\rVert_F.
\]

One of the following statements must be true:
\[
\lVert X\tilde{Q}^{\top}\tilde{Q}X-Y\tilde{Q}^{\top}\tilde{Q}Y\rVert_F\leq \lVert X\tilde{Q}^{\top}\tilde{Q}X + X\tilde{Q}^{\top}\tilde{Q}Y - Y\tilde{Q}^{\top}\tilde{Q}X -Y\tilde{Q}^{\top}\tilde{Q}Y \rVert_F
\]
or
\[
\lVert X\tilde{Q}^{\top}\tilde{Q}X-Y\tilde{Q}^{\top}\tilde{Q}Y\rVert_F\leq \lVert X\tilde{Q}^{\top}\tilde{Q}X - X\tilde{Q}^{\top}\tilde{Q}Y + Y\tilde{Q}^{\top}\tilde{Q}X -Y\tilde{Q}^{\top}\tilde{Q}Y \rVert_F.
\]
Without loss of generality we only treat the first case. Thus we obtain
\[
\lVert X\tilde{Q}^{\top}\tilde{Q}X-Y\tilde{Q}^{\top}\tilde{Q}Y\rVert_F\leq \lVert (X-Y)\tilde{Q}^{\top}\tilde{Q}(X+Y) \rVert_F.
\]

The matrix $\tilde{Q}^{\top}\tilde{Q}(\tilde{Q}^{\top}\tilde{Q})^{\top}$ is symmetric and thus its spectral decomposition can be written as $O\Lambda O^{\top}$ for an orthogonal matrix $O$. Let $\lambda_{\text{max}}$ be the largest eigenvalue of $\tilde{Q}^{\top}\tilde{Q}(\tilde{Q}^{\top}\tilde{Q})^{\top}$. 
Since $\tilde{Q}^{\top}\tilde{Q}(\tilde{Q}^{\top}\tilde{Q})^{\top}$ has only non-negative eigenvalues, $\lambda_{\text{max}}I-\Lambda$ is positive definite and so is $(X-Y)(X+Y)O(\lambda_{\text{max}}I-\Lambda)O^{\top}(X+Y)^{\top}(X-Y)^{\top}$. Since the trace of a matrix is the sum of its eigenvalues we have
\begin{equation*}
    \begin{aligned}
     & \operatorname{Tr}[(X-Y)(X+Y)O(\lambda_{\text{max}}I-\Lambda)O^{\top}(X+Y)^{\top}(X-Y)^{\top}]\geq 0
    \end{aligned}
\end{equation*}

and thus 
\[
\lambda_{\text{max}}\operatorname{Tr}[(X-Y)(X+Y)(X+Y)^{\top}(X-Y)^{\top}] \geq \operatorname{Tr}[(X-Y)(X+Y)O\Lambda O^{\top}(X+Y)^{\top}(X-Y)^{\top}]
\]

Using the above facts, we obtain

\begin{equation*}
    \begin{aligned}
     \lVert (X-Y)\tilde{Q}^{\top}\tilde{Q}(X+Y) \rVert_F^2{}
     & = \operatorname{Tr}[(X-Y)\tilde{Q}^{\top}\tilde{Q}(X+Y)(X+Y)^{\top}(\tilde{Q}^{\top}\tilde{Q})^{\top}(X-Y)^{\top}]\\
     & = \operatorname{Tr}[(X-Y)(X+Y)\tilde{Q}^{\top}\tilde{Q}(\tilde{Q}^{\top}\tilde{Q})^{\top}(X+Y)^{\top}(X-Y)^{\top}]\\
     & \leq \lambda_{\text{max}}\operatorname{Tr}[(X-Y)(X+Y)(X+Y)^{\top}(X-Y)^{\top}]\\
     & = \lambda_{\text{max}}\lVert (X-Y)(X+Y) \rVert_F^2.
    \end{aligned}
\end{equation*}

Hence we have

\begin{equation*}
    \begin{aligned}
     \abs{\vecop(X\tilde{Q}^{\top}\tilde{Q}X-Y\tilde{Q}^{\top}\tilde{Q}Y)} {}
     & \leq \sqrt{\lambda_{\text{max}}} \abs{\vecop((X-Y)(X+Y))} \\
     & = \sqrt{\lambda_{\text{max}}} \abs{(I\otimes(X-Y)\vecop(X+Y))}\\
     &\leq \sqrt{\lambda_{\text{max}}} \lVert I\otimes (X-Y) \rVert_F \abs{\vecop(X+Y)}\\
     & = \sqrt{\lambda_{\text{max}}} \sqrt{d} \lVert (X-Y) \rVert_F \abs{\vecop(X+Y)} \\
     & \leq \sqrt{\lambda_{\text{max}}} \sqrt{d} \abs{\vecop(X)-\vecop(Y)} (\abs{\vecop(X)}+\abs{\vecop(Y)})
    \end{aligned}
\end{equation*}

leading to $C_2=2\sqrt{d}\sqrt{\lambda_{\text{max}}} $.
Now an application of \cite[Theorem B2]{AJ19} yields the desired result.  \qed\\

\appendix

\section{Detailed proof of the main result}\label{A:main_result}
\emph{Proof:}
In order to prove that $\pi^*$ is indeed the optimal portfolio strategy, we show that for 
\[
G(x_0,\Sigma_0):=\frac{x_0^\gamma}{\gamma}\expo\Big(\int_0^T \gamma r_s + \operatorname{Tr}[f(\psi)(T-s)\Sigma_0 +\psi(T-s)NN^{\top}]ds\Big), 
\]
we have
\begin{enumerate}
    \item $\EX^{x_0,\Sigma_0}{\big[\frac{1}{\gamma}(X_T^{\pi^*})^{\gamma}\big]}=G(x_0,\Sigma_0)$ for $\pi_t^*=\frac{1}{1-\gamma}(v+2\psi(T-t)Q^{\top}\rho)$,\label{a1A}
    \item $\EX^{x_0,\Sigma_0}{\big[\frac{1}{\gamma}(X_T^{\pi})^{\gamma}\big]}\leq G(x_0,\Sigma_0)$ for every other admissible strategy.\label{b1A}
\end{enumerate}

We start with equation \eqref{a1A}.
The SDE for the wealth process can be solved explicitly and we obtain for an arbitrary admissible portfolio strategy $\pi_t$
\[
X_T^{\pi}=X_0^{\pi}\expo\Big(\int_0^T(r_s+\pi_s^{\top}\Sigma_s v-\frac{1}{2}\lVert \pi_s^{\top}\Sigma_s^{1/2} \rVert_{2}^2 )ds + \int_0^T \pi_s^{\top}\Sigma_s^{1/2}dW_s^S\Big)
\]
with $X_0^\pi=x_0$.

Since $\pi_t^*$ is continuous by the continuity of $\psi$, it is also bounded and thus we can define a new probability measure ${\Q}$ with Radon-Nikodym density  
\begin{equation}\label{density QA}
{Z}_0:=\frac{d{\Q}}{d\mathbb{P}}|_{\mathcal{F}_0}=\expo\Big(\gamma\int_0^T (\pi_s^{*})^{\top}\Sigma_s^{1/2}dW_{s}^S-\frac{\gamma^2}{2}\int_0^T \lVert (\pi_s^{*})^{\top}\Sigma_s^{1/2} \rVert_{2}^2 ds\Big),
\end{equation}
 by Lemma \ref{A2}.
Then we obtain
\begin{equation*}
\begin{aligned}
x_0^{-\gamma}\EX_{x_0,\Sigma_0}[{(X_T^{\pi^*})^{\gamma}}]{}
    &=\EX_{x_0,\Sigma_0}\Big[\expo\Big(\gamma\int_0^T(r_s+{(\pi_s^{*})^{\top}}\Sigma_s v-\frac{1}{2}\lVert (\pi_s^{*})^{\top}\Sigma_s^{1/2} \rVert_{2}^2 )ds\\
    &\qquad\qquad\quad+ \gamma\int_0^T (\pi_s^{*})^{\top}\Sigma_s^{1/2}dW_s^S \Big)\Big]\\
    & =\EX_{x_0,\Sigma_0}\Big[\expo\Big(\gamma\int_0^T(r_s+(\pi_s^{*})^{\top}\Sigma_s v-\frac{1}{2}\lVert (\pi_s^{*})\Sigma_s^{1/2} \rVert_{2}^2 )ds \\
    &\qquad\qquad\quad +\frac{\gamma^2}{2}\int_0^T \lVert (\pi_s^{*})^{\top}\Sigma_s^{1/2} \rVert_{2}^2 ds\Big)\\
    & \qquad\qquad\times \expo\Big(\gamma\int_0^T (\pi_s^{*})^{\top}\Sigma_s^{1/2}dW_s^S-\frac{\gamma^2}{2}\int_0^T \lVert (\pi_s^{*})^{\top}\Sigma_s^{1/2} \rVert_{2}^2 ds\Big)\Big]\\
    & =\EX_{x_0,\Sigma_0}^{\Q}\Big[\expo\Big(\gamma\int_0^T(r_s+(\pi_s^{*})^{\top}\Sigma_s v-\frac{1}{2}\lVert (\pi_s^{*})^{\top}\Sigma_s^{1/2}\rVert_{2}^2 )ds\\
    & \qquad\qquad\quad+\frac{\gamma^2}{2}\int_0^T \lVert (\pi_s^{*})^{\top}\Sigma_s^{1/2} \rVert_{2}^2 ds\Big)\\
    & =\EX_{x_0,\Sigma_0}^{\Q}\Big[\expo\Big(\gamma\int_0^T\big[r_s+(\pi_s^{*})^{\top}\Sigma_s v
    +\frac{\gamma-1}{2} (\pi_s^{*})^{\top}\Sigma_s\pi_s^* \big]ds\Big)\Big]\\
    & =\EX_{x_0,\Sigma_0}^{\Q}\Big[\expo\Big(\int_0^T \gamma r_s ds+\int_0^T \operatorname{Tr}\big[\big(\gamma v(\pi_s^{*})^{\top}
    +\frac{\gamma(\gamma-1)}{2} \pi_s^*(\pi_s^{*})^{\top} \big)\Sigma_s\big] ds\Big)\Big].\\
\end{aligned}
\end{equation*}
In the following, we denote the matrix-valued deterministic process $F_s$ by 
\[
F_s:=\gamma v(\pi_s^{*})^{\top}
    +\frac{\gamma(\gamma-1)}{2} \pi_s^*(\pi_s^{*})^{\top}.
\]
Inserting the optimal strategy $\pi_t^*=\frac{1}{1-\gamma}(v+2\psi(T-t)Q^{\top}\rho)$, we obtain
\begin{equation*}
    \begin{aligned}
    F_s &
    = \frac{\gamma}{1-\gamma}(\frac{1}{2}vv^{\top}+v\rho^{\top}Q\psi(T-s)-\psi(T-s)Q^{\top}\rho v^{\top}-2\psi(T-s)Q^{\top}\rho\rho^{\top}Q\psi(T-s))\\
    & = f(\psi)(T-s)-\psi(T-s)M-M^{\top}\psi(T-s)-2\psi(T-s)Q^{\top}Q\psi(t-s)\\
    &\quad + \frac{\gamma}{\gamma-1}(2\psi(T-s)Q^{\top}\rho v^{\top}+4\psi(T-s)Q^{\top}\rho\rho^{\top}Q\psi(T-s)).
    \end{aligned}
\end{equation*}

Under the probability measure ${\Q}$ defined in \eqref{density Q}, the process 
\[
\tilde{W}_{t}^{\sigma}= W_t^{\sigma}-\gamma\Sigma_t^{1/2}\pi_t^{*}\rho^{\top}
\]
is a $d\times d$-dimensional brownian motion by Girsanov's theorem. 
The dynamics of the volatility process $\Sigma$ under $\hat{\Q}$ can thus be written as

\begin{equation}\label{Sigma_subst}
    \begin{aligned}
    \Sigma_t&
    =\Sigma_0+\int_0^t K(t-s)(NN^{\top}+M\Sigma_s+\Sigma_s M^{\top}+\frac{\gamma}{1-\gamma}\Sigma_s v\rho^{\top}Q\\
    &\quad +\frac{2\gamma}{1-\gamma}\Sigma_s\psi(T-s)Q^{\top}\rho\rho^{\top}Q+\frac{\gamma}{1-\gamma}Q^{\top}\rho v^{\top}\Sigma_s+\frac{2\gamma}{1-\gamma}Q^{\top}\rho\rho^{\top}Q\psi(T-s)\Sigma_s)ds\\
    &\quad + \int_0^t K(t-s)\Sigma_s^{1/2}d\tilde{W}_s^{\sigma}Q+K(t-s)Q^{\top}(d\tilde{W}_s^{\sigma})^{\top}\Sigma_s^{1/2}.
    \end{aligned}
\end{equation}

Let $L$ be the resolvent of the first kind of the integral kernel $K$.
Then by the associativity of the convolution and applying the fundamental theorem of calculus we obtain
\[
(\psi\ast L)(t)=((f(\psi)\ast K)\ast L)(t)=(f(\psi)\ast(K\ast L))(t)=(f(\psi)\ast I)(t)=\int_0^t f(\psi)(s) ds
\]
The fundamental theorem of calculus thus yields  
\[
f(\psi)(t)=(\psi\ast L)'(t).
\]
Thus we can write $F_s$ as
\begin{equation*}
    \begin{aligned}
    F_s{}
    & = (\psi\ast L)'(T-s)-\psi(T-s)M-M^{\top}\psi(T-s)-2\psi(T-s)Q^{\top}Q\psi(t-s)\\
    &\quad + \frac{\gamma}{\gamma-1}(2\psi(T-s)Q^{\top}\rho v^{\top}+4\psi(T-s)Q^{\top}\rho\rho^{\top}Q\psi(T-s)),
    \end{aligned}
\end{equation*}
and consequently we have

\begin{equation*}
    \begin{aligned}
    \int_0^T &\operatorname{Tr}[F_s\Sigma_s]ds{}\\
    & = \operatorname{Tr}[\int_0^T (\psi\ast L)'(T-s)\Sigma_s ds]-\operatorname{Tr}[\int_0^T\psi(T-s)M\Sigma_sds]\\
    &\quad-\operatorname{Tr}[\int_0^T M^{\top}\psi(T-s)\Sigma_s ds]-\operatorname{Tr}[\int_0^T 2\psi(T-s)Q^{\top}Q\psi(t-s)\Sigma_s ds]\\
    &\quad - \operatorname{Tr}[\int_0^T\frac{2\gamma}{1-\gamma} \psi(T-s)Q^{\top}\rho v^{\top}\Sigma_s ds]-\operatorname{Tr}[\int_0^T\frac{4\gamma}{1-\gamma}\psi(T-s)Q^{\top}\rho\rho^{\top}Q\psi(T-s)\Sigma_sds].
    \end{aligned}
\end{equation*}

We consider the term $\operatorname{Tr}[\int_0^T (\psi\ast L)'(T-s)\Sigma_s ds]$ {and substitute \eqref{Sigma_subst} for $\Sigma$}.
This yields

\begin{equation*}
\begin{aligned}
 {}
    & \operatorname{Tr}[\int_0^T (\psi\ast L)'(T-s)\Sigma_s ds] =  \operatorname{Tr}[\left((\psi\ast L)'\ast \Sigma\right)(T)]\\
    & = \operatorname{Tr}[\Big((\psi\ast L)'\ast \big(\Sigma_0+ K\ast\big(NN^{\top}+M\Sigma+\Sigma M^{\top}+\frac{\gamma}{1-\gamma}\Sigma v\rho^{\top}Q+\frac{2\gamma}{1-\gamma}\Sigma\psi^{T-}Q^{\top}\rho\rho^{\top}Q\\
    & \quad +\frac{\gamma}{1-\gamma}Q^{\top}\rho v^{\top}\Sigma+\frac{2\gamma}{1-\gamma}Q^{\top}\rho\rho^{\top}Q\psi^{T-}\Sigma+ \Sigma_s^{1/2}d\tilde{W}_s^{\sigma}Q+Q^{\top}(d\tilde{W}_s^{\sigma})^{\top}\Sigma_s^{1/2}\big)\big)\Big)(T)]\\
    &=  \operatorname{Tr}[\left((\psi\ast L)'\ast\Sigma_0\right)(T)]\qquad (\textbf{I})\\ 
    & \quad +  \operatorname{Tr}[\Big((\psi\ast L)'\ast  \big(K\ast\big(NN^{\top}+M\Sigma+\Sigma M^{\top}+\frac{\gamma}{1-\gamma}\Sigma v\rho^{\top}Q\\
    & \qquad +\frac{2\gamma}{1-\gamma}\Sigma\psi^{T-}Q^{\top}\rho\rho^{\top}Q+\frac{\gamma}{1-\gamma}Q^{\top}\rho v^{\top}\Sigma+\frac{2\gamma}{1-\gamma}Q^{\top}\rho\rho^{\top}Q\psi^{T-}\Sigma \big)\big)\Big)(T)]\quad  (\textbf{II})\\
    & \quad + \operatorname{Tr}[\Big((\psi\ast L)'\ast \big(K\ast ( \Sigma_s^{1/2}d\tilde{W}_s^{\sigma}Q+Q^{\top}(d\tilde{W}_s^{\sigma})^{\top}\Sigma_s^{1/2})\big)\Big)(T)]. \qquad (\textbf{III})
\end{aligned}
\end{equation*}

{Here $\psi^{T-}(s):= \psi(T-s)$.}
We now simplify the terms $(\textbf{I})$, $(\textbf{II})$ and $(\textbf{III})$ using stochastic calculus of convolutions and resolvents.

ad (\textbf{I}):
\begin{equation*}
\begin{aligned}
\operatorname{Tr}\Big[\Big((\psi\ast L)'\ast \Sigma_0\Big)(T)\Big] {}
    &  = \operatorname{Tr}\Big[\int_0^T (\psi\ast L)'(T-s) \Sigma_0 ds\Big] \\
    & = \operatorname{Tr}\Big[\int_0^T f(\psi)(T-s)\Sigma_0 ds\Big]
\end{aligned}
\end{equation*}

ad (\textbf{II}):
\allowdisplaybreaks

\begin{align*}
{}
    & \operatorname{Tr}\Big[ \Big((\psi\ast L)' \ast  \big(K\ast\big(NN^{\top}+M\Sigma+\Sigma M^{\top}+\frac{\gamma}{1-\gamma}\Sigma v\rho^{\top}Q\\
    & \quad   +\frac{2\gamma}{1-\gamma}\Sigma\psi^{T-}Q^{\top}\rho\rho^{\top}Q  +\frac{\gamma}{1-\gamma}Q^{\top}\rho v^{\top}\Sigma+\frac{2\gamma}{1-\gamma}Q^{\top}\rho\rho^{\top}Q\psi^{T-}\Sigma\big)\big)\Big)(T)\Big]\\
    &= \operatorname{Tr}\Big[ \Big(((\psi\ast L)\ast K)'\ast  \big(NN^{\top}+M\Sigma+\Sigma M^{\top}+\frac{\gamma}{1-\gamma}\Sigma v\rho^{\top}Q\\
    & \quad  +\frac{2\gamma}{1-\gamma}\Sigma\psi^{T-}Q^{\top}\rho\rho^{\top}Q  +\frac{\gamma}{1-\gamma}Q^{\top}\rho v^{\top}\Sigma+\frac{2\gamma}{1-\gamma}Q^{\top}\rho\rho^{\top}Q\psi^{T-}\Sigma\big)\Big)(T)\Big]\\
    &= \operatorname{Tr}\Big[ \Big((\psi\ast (L\ast K))'\ast  \big(NN^{\top}+M\Sigma+\Sigma M^{\top}+\frac{\gamma}{1-\gamma}\Sigma v\rho^{\top}Q\\
    & \quad  +\frac{2\gamma}{1-\gamma}\Sigma\psi^{T-}Q^{\top}\rho\rho^{\top}Q  +\frac{\gamma}{1-\gamma}Q^{\top}\rho v^{\top}\Sigma+\frac{2\gamma}{1-\gamma}Q^{\top}\rho\rho^{\top}Q\psi^{T-}\Sigma\big)\Big)(T)\Big]\\
    &= \operatorname{Tr}\Big[ \Big((\psi\ast I)'\ast  \big(NN^{\top}+M\Sigma+\Sigma M^{\top}+\frac{\gamma}{1-\gamma}\Sigma v\rho^{\top}Q\\
    & \quad +\frac{2\gamma}{1-\gamma}\Sigma\psi^{T-}Q^{\top}\rho\rho^{\top}Q +\frac{\gamma}{1-\gamma}Q^{\top}\rho v^{\top}\Sigma+\frac{2\gamma}{1-\gamma}Q^{\top}\rho\rho^{\top}Q\psi^{T-}\Sigma\big)\Big)(T)\Big]\\
    &= \operatorname{Tr}\Big[ \Big(\psi\ast\big(NN^{\top}+M\Sigma+\Sigma M^{\top}+\frac{\gamma}{1-\gamma}\Sigma v\rho^{\top}Q+\frac{2\gamma}{1-\gamma}\Sigma\psi^{T-}Q^{\top}\rho\rho^{\top}Q\\
    & \quad +\frac{\gamma}{1-\gamma}Q^{\top}\rho v^{\top}\Sigma+\frac{2\gamma}{1-\gamma}Q^{\top}\rho\rho^{\top}Q\psi^{T-}\Sigma\big)\Big)(T)\Big]\\
    & = \operatorname{Tr}[\int_0^T \psi(T-s)NN^{\top} ds]+\operatorname{Tr}[\int_0^T\psi(T-s)M\Sigma_sds]
    +\operatorname{Tr}[\int_0^T M^{\top}\psi(T-s)\Sigma_s ds]\\
    &\quad + \operatorname{Tr}[\int_0^T\frac{2\gamma}{1-\gamma} \psi(T-s)Q^{\top}\rho v^{\top}\Sigma_s ds]+\operatorname{Tr}[\int_0^T\frac{4\gamma}{1-\gamma}\psi(T-s)Q^{\top}\rho\rho^{\top}Q\psi(T-s)\Sigma_sds].
\end{align*}

ad (\textbf{III}):

The processes $dM_1:=\Sigma_s^{1/2}Q d\tilde{W}_s^{\sigma}$ and $dM_2:=Q^{\top}\Sigma_s^{1/2}(d\tilde{W}_s^{\sigma})^{\top}$ are both continuous local martingales and $\langle M_1 \rangle_t$ and $\langle M_2 \rangle_t$ are locally bounded. Thus we can apply Lemma \ref{lem:assoc} to obtain
    
\begin{equation*}
\begin{aligned}
{}
    & \operatorname{Tr}\Big[\Big((\psi\ast L)'\ast \big(K\ast ( \Sigma_s^{1/2}d\tilde{W}_s^{\sigma}Q+Q^{\top}(d\tilde{W}_s^{\sigma})^{\top}\Sigma_s^{1/2})\big)\Big)(T)\Big]\\
    &= \operatorname{Tr}\Big[\Big((\psi\ast (L\ast K))'\ast ( \Sigma_s^{1/2}d\tilde{W}_s^{\sigma}Q+Q^{\top}(d\tilde{W}_s^{\sigma})^{\top}\Sigma_s^{1/2})\Big)(T)\Big]\\
    &= \operatorname{Tr}\Big[\Big(\psi\ast ( \Sigma_s^{1/2}d\tilde{W}_s^{\sigma}Q+Q^{\top}(d\tilde{W}_s^{\sigma})^{\top}\Sigma_s^{1/2})\Big)(T)\Big]\\
    &= \operatorname{Tr}\Big[\int_0^T \psi(T-s)( \Sigma_s^{1/2}d\tilde{W}_s^{\sigma}Q+Q^{\top}(d\tilde{W}_s^{\sigma})^{\top}\Sigma_s^{1/2})\Big].
\end{aligned}
\end{equation*}

Combining the above results, we end up with

\begin{equation*}
    \begin{aligned}
    \int_0^T \operatorname{Tr}[F_sV_s]ds{}
    &= \operatorname{Tr}[\int_0^T f(\psi)(T-s)\Sigma_0 ds]+\operatorname{Tr}[\int_0^T \psi(T-s)NN^{\top} ds]\\ &\qquad-\operatorname{Tr}[\int_0^T 2\psi(T-s)Q^{\top}Q\psi(t-s)\Sigma_s ds]\\ 
    &\qquad +\operatorname{Tr}[\int_0^T \psi(T-s)( \Sigma_s^{1/2}d\tilde{W}_s^{\sigma}Q+Q^{\top}(d\tilde{W}_s^{\sigma})^{\top}\Sigma_s^{1/2})].
    \end{aligned}
\end{equation*}

Hence we obtain

\begin{equation*}
    \begin{aligned}
    {}
    & x_0^{-\gamma}\EX_{x_0,\Sigma_0}{(X_T^{\pi^*})^{\gamma}}=\EX_{x_0,\Sigma_0}^{\Q}\Big[\expo\Big(\int_0^T \gamma r_s ds+\int_0^T \operatorname{Tr}[F_s \Sigma_s] ds\Big)\Big]\\
    & = \EX_{x_0,\Sigma_0}^{\Q}\Big[\expo\Big(\int_0^T \gamma r_s ds + \operatorname{Tr}[\int_0^T f(\psi)(T-s)\Sigma_0 ds]+\operatorname{Tr}[\int_0^T \psi(T-s)NN^{\top} ds] \\
    &\qquad -\operatorname{Tr}[\int_0^T 2\psi(T-s)Q^{\top}Q\psi(t-s)\Sigma_s ds]\\
    &\qquad +\operatorname{Tr}[\int_0^T \psi(T-s) \Sigma_s^{1/2}d\tilde{W}_s^{\sigma}Q]+\operatorname{Tr}[\int_0^T \psi(T-s)Q^{\top}(d\tilde{W}_s^{\sigma})^{\top}\Sigma_s^{1/2}]\Big)\Big]\\
    &= \expo\Big(\int_0^T \gamma r_s ds + \operatorname{Tr}[\int_0^T f(\psi)(T-s)\Sigma_0 ds]+\operatorname{Tr}[\int_0^T \psi(T-s)NN^{\top} ds]\Big)\\
    &\quad \times\EX_{x_0,\Sigma_0}^{\Q}\Big[\expo\Big(-2\operatorname{Tr}[\int_0^T Q\psi(T-s)\Sigma_s\psi(t-s)Q^{\top} ds]  + 2\operatorname{Tr}[\int_0^T Q\psi(T-s)\Sigma_s^{1/2}d\tilde{W}_s^{\sigma}] \Big)\Big].\\
    \end{aligned}
\end{equation*}

Since $\psi$ is continuous, it is bounded and therefore the stochastic exponential is a true $\Q$- martingale with expectation $1$ by Lemma \ref{A1}. Thus we get
\[
\EX_{x_0,\Sigma_0}\Big[{\frac{1}{\gamma}(X_T^{\pi^*})^{\gamma}}\Big]=\frac{x_0^\gamma}{\gamma}\expo\Big(\int_0^T \gamma r_s + \operatorname{Tr}[f(\psi)(T-s)\Sigma_0 +\psi(T-s)NN^{\top}]ds\Big).
\]
This completes the first part of the proof.

It remains to show the inequality \eqref{b1A} for arbitrary admissible portfolio strategies.

To this end, let $\pi_t$ be an arbitrary admissible portfolio strategy. Define $\hat{\pi}:=\pi-\pi^*$ and write $\pi_t=\pi_t^* +\hat{\pi}_t$, where $\pi_t^*=\frac{1}{1-\gamma}(v+2\psi(T-t)Q^{\top}\rho$. Since $\pi_t$ is bounded by assumption, we can define a new probability measure $\hat{\Q}$ with Radon-Nikodym density  
\begin{equation}\label{density hatQA}
\hat{Z}_0:=\frac{d\hat{\Q}}{d\mathbb{P}}|_{\mathcal{F}_0}=\expo\Big(\gamma\int_0^T \pi_s^{\top}\Sigma_s^{1/2}dW_s^{S}-\frac{\gamma^2}{2}\int_0^T \lVert \pi_s^{\top}\Sigma_s^{1/2} \rVert_{2}^2 ds\Big)
\end{equation}
by Lemma \ref{A2}.   

Analogously to the first part we obtain

\begin{equation*}
    \begin{aligned}
    & x_0^{-\gamma}\EX_{x_0,\Sigma_0}{(X_T^{\pi})^{\gamma}}=\EX_{x_0,\Sigma_0}^{\hat{\Q}}\Big[\expo\Big(\int_0^T \gamma r_s ds+\int_0^T \operatorname{Tr}\big[\big(\gamma v\pi_s^{\top}
    +\frac{\gamma(\gamma-1)}{2} \pi_s\pi_s^{\top} \big)\Sigma_s\big] ds\Big)\Big].
    \end{aligned}
\end{equation*}

We define 
\begin{equation*}
\begin{aligned}
\hat{F}_s {}
    &:= \gamma v\pi_s^{\top}
    +\frac{\gamma(\gamma-1)}{2} \pi_s\pi_s^{\top} = \gamma v ((\pi_s^*)^{\top}+\hat{\pi}_s^{\top}) + \frac{\gamma(\gamma-1)}{2}(\pi_s^*+\hat{\pi}_s)((\pi_s^*)^{\top}+\hat{\pi}_s^{\top})\\
    & = \gamma v(\pi_s^*)^{\top}+\frac{\gamma(\gamma-1)}{2}\pi_s^*(\pi_s^*)^{\top} + \gamma v \hat{\pi}_s^{\top} + \gamma(\gamma-1)\pi_s^*\hat{\pi}_s^{\top} + \frac{\gamma(\gamma-1)}{2}\hat{\pi}_s\hat{\pi}_s^{\top}\\
    &= F_s + \gamma v \hat{\pi}_s^{\top} + \gamma(\gamma-1)\pi_s^*\hat{\pi}_s^{\top} + \frac{\gamma(\gamma-1)}{2}\hat{\pi}_s\hat{\pi}_s^{\top},
\end{aligned}
\end{equation*}

with

\begin{equation*}
    \begin{aligned}
    F_s{}
    & = (\psi\ast L)'(T-s)-\psi(T-s)M-M^{\top}\psi(T-s)-2\psi(T-s)Q^{\top}Q\psi(t-s)\\
    &\quad + \frac{\gamma}{\gamma-1}(2\psi(T-s)Q^{\top}\rho v^{\top}+4\psi(T-s)Q^{\top}\rho\rho^{\top}Q\psi(T-s)).
    \end{aligned}
\end{equation*}

Under the probability measure $\hat{\Q}$ defined in \eqref{density hatQ}, the process 
\[
\hat{W}_{t}^{\sigma}=W_t^{\sigma}-\gamma\Sigma_t^{1/2}\pi_t\rho^{\top}= W_t^{\sigma}-\gamma\Sigma_t^{1/2}\pi_t^{*}\rho^{\top}-\gamma\Sigma_t^{1/2}\hat{\pi}_t\rho^{\top}
\]
is a $d\times d$-dimensional brownian motion by Girsanov's theorem. 

Then the dynamics of the volatility process $\Sigma$ under $\hat{\Q}$ can be written as

\begin{equation*}
    \begin{aligned}
    \Sigma_t&
    =\Sigma_0+\int_0^t K(t-s)(NN^{\top}+M\Sigma_s+\Sigma_s M^{\top}+\frac{\gamma}{1-\gamma}\Sigma_s v\rho^{\top}Q\\
    &\quad +\frac{2\gamma}{1-\gamma}\Sigma_s\psi(T-s)Q^{\top}\rho\rho^{\top}Q+\frac{\gamma}{1-\gamma}Q^{\top}\rho v^{\top}\Sigma_s+\frac{2\gamma}{1-\gamma}Q^{\top}\rho\rho^{\top}Q\psi(T-s)\Sigma_s\\
    &\quad + \gamma\Sigma_s\hat{\pi}\rho^{\top}Q+\gamma Q^{\top}\rho\hat{\pi}^{\top}\Sigma_s) + \int_0^t K(t-s)\Sigma_s^{1/2}d\hat{W}_s^{\sigma}Q+Q^{\top}(d\hat{W}_s^{\sigma})^{\top}\Sigma_s^{1/2}K(t-s)
    \end{aligned}
\end{equation*}

Therefore, we get

\begin{equation*}
\begin{aligned}
 {}
    &\operatorname{Tr}[\int_0^T (\psi\ast L)'(T-s)\Sigma_s ds] \\
    &= \operatorname{Tr}[\int_0^T f(\psi)(T-s)\Sigma_0 ds]+\operatorname{Tr}[\int_0^T \psi(T-s)NN^{\top} ds]\\
    &\quad+\operatorname{Tr}[\int_0^T\psi(T-s)M\Sigma_sds]
    +\operatorname{Tr}[\int_0^T M^{\top}\psi(T-s)\Sigma_s ds]\\
    &\quad + \operatorname{Tr}[\int_0^T\frac{2\gamma}{1-\gamma} \psi(T-s)Q^{\top}\rho v^{\top}\Sigma_s ds]+\operatorname{Tr}[\int_0^T\frac{4\gamma}{1-\gamma}\psi(T-s)Q^{\top}\rho\rho^{\top}Q\psi(T-s)\Sigma_sds]\\
    &\quad + \operatorname{Tr}[\int_0^T \psi(T-s)\gamma\Sigma_s\hat{\pi}\rho^{\top}Qds]+\operatorname{Tr}[\int_0^T \psi(T-s)\gamma Q^{\top}\rho\hat{\pi}^{\top}\Sigma_s ds ]\\
    &\quad +\operatorname{Tr}[\int_0^T \psi(T-s) \Sigma_s^{1/2}d\hat{W}_s^{\sigma}Q+\int_0^T \psi(T-s)Q^{\top}(d\hat{W}_s^{\sigma})^{\top}\Sigma_s^{1/2}].
\end{aligned}
\end{equation*}

Combining the above results, we obtain

\begin{equation*}
\begin{aligned}
&\int_0^T \operatorname{Tr}[\hat{F}_s \Sigma_s] ds {}\\
    &\quad= \operatorname{Tr}[\int_0^T [F_s + \gamma v \hat{\pi}_s^{\top} + \gamma(\gamma-1)\pi_s^*\hat{\pi}_s^{\top} + \frac{\gamma(\gamma-1)}{2}\hat{\pi}_s\hat{\pi}_s^{\top}]\Sigma_s ds]\\
    &\quad= \operatorname{Tr}[\int_0^T f(\psi)(T-s)\Sigma_0 ds]+\operatorname{Tr}[\int_0^T \psi(T-s)NN^{\top} ds]\\
    &\qquad -\operatorname{Tr}[\int_0^T 2\psi(T-s)Q^{\top}Q\psi(t-s)\Sigma_s ds] \\
    &\qquad +\operatorname{Tr}[\int_0^T \psi(T-s)( \Sigma_s^{1/2}d\hat{W}_s^{\sigma}Q+Q^{\top}(d\hat{W}_s^{\sigma})^{\top}\Sigma_s^{1/2})]+2\operatorname{Tr}[\int_0^T \psi(T-s)\gamma Q^{\top}\rho\hat{\pi}^{\top}\Sigma_s ds ]\\
    &\qquad+  \operatorname{Tr}[\int_0^T \gamma v \hat{\pi}_s^{\top}\Sigma_s ds] + \operatorname{Tr}[\int_0^T\gamma(\gamma-1)\pi_s^*\hat{\pi}_s^{\top}\Sigma_s ds] + \operatorname{Tr}[\int_0^T\frac{\gamma(\gamma-1)}{2}\hat{\pi}_s\hat{\pi}_s^{\top}\Sigma_s ds]
\end{aligned}
\end{equation*}

Since

\begin{equation*}
    \begin{aligned}
     \operatorname{Tr}[\int_0^T \gamma v \hat{\pi}_s^{\top}\Sigma_s ds]  + 2\operatorname{Tr}[\int_0^T \psi(T-s)\gamma Q^{\top}\rho\hat{\pi}^{\top}\Sigma_s ds ]
    &= \operatorname{Tr}[\int_0^T \gamma(v+2\psi(T-s)Q^{\top}\rho)\hat{\pi}_s^{\top}\Sigma_s ds]\\
    &= \operatorname{Tr}[\int_0^T\gamma(1-\gamma)\pi_s^*\hat{\pi}_s^{\top}\Sigma_s ds],
    \end{aligned}
\end{equation*}

three terms in the above sum cancel out and we end up with

\begin{equation*}
\begin{aligned}
&\int_0^T \operatorname{Tr}[\hat{F}_s \Sigma_s] ds {}\\
    &\quad= \operatorname{Tr}[\int_0^T f(\psi)(T-s)\Sigma_0 ds]+\operatorname{Tr}[\int_0^T \psi(T-s)NN^{\top} ds]\\
    &\qquad-\operatorname{Tr}[\int_0^T 2\psi(T-s)Q^{\top}Q\psi(t-s)\Sigma_s ds] \\
    &\qquad +\operatorname{Tr}[\int_0^T \psi(T-s)( \Sigma_s^{1/2}d\hat{W}_s^{\sigma}Q+Q^{\top}(d\hat{W}_s^{\sigma})^{\top}\Sigma_s^{1/2})]+ \operatorname{Tr}[\int_0^T\frac{\gamma(\gamma-1)}{2}\hat{\pi}_s\hat{\pi}_s^{\top}\Sigma_s ds]
\end{aligned}
\end{equation*}

Thus, for the expectation we get 

\begin{equation*}
    \begin{aligned}
    {}
    & x_0^{-\gamma}\EX_{x_0,\Sigma_0}[{(X_T^{\pi})^{\gamma}}]=\EX_{x_0,\Sigma_0}^{\hat{\Q}}\Big[\expo\Big(\int_0^T \gamma r_s ds+\int_0^T \operatorname{Tr}[\hat{F}_s\Sigma_s] ds\Big)\Big]\\
    &=\expo\Big(\int_0^T \gamma r_s ds + \operatorname{Tr}[\int_0^T f(\psi)(T-s)\Sigma_0 ds]+\operatorname{Tr}[\int_0^T \psi(T-s)NN^{\top} ds] \Big)\\
    &\qquad\times\EX_{x_0,\Sigma_0}^{\hat{\Q}}\Big[\expo\Big(-\operatorname{Tr}[\int_0^T 2\psi(T-s)Q^{\top}Q\psi(t-s)\Sigma_s ds]  + 2\operatorname{Tr}[\int_0^T Q\psi(T-s)\Sigma_s^{1/2}d\hat{W}_s^{\sigma}] \\
    & \qquad\qquad\qquad + \operatorname{Tr}[\int_0^T\frac{\gamma(\gamma-1)}{2}\hat{\pi}_s\Sigma_s\hat{\pi}_s^{\top} ds]\Big)\Big]\\
    &\leq \expo\Big(\int_0^T \gamma r_s ds + \operatorname{Tr}[\int_0^T f(\psi)(T-s)\Sigma_0 ds]+\operatorname{Tr}[\int_0^T \psi(T-s)NN^{\top} ds] \Big)\\ 
    &\quad \times\EX_{x_0,\Sigma_0}^{\hat{\Q}}\Big[\expo\Big(-2\operatorname{Tr}[\int_0^T Q\psi(T-s)\Sigma_s\psi(t-s)Q^{\top} ds]  + 2\operatorname{Tr}[\int_0^T Q\psi(T-s)\Sigma_s^{1/2}d\hat{W}_s^{\sigma}] \Big)\Big].
    \end{aligned}
\end{equation*}

The last inequality follows from the fact that $\Sigma_s$ is positive definite and $\gamma\in (0,1)$. Since the stochastic exponential is a $\hat{\Q}$-martingale with expectation $1$ by Lemma \ref{A1}, we finally obtain

\[
\EX_{x_0,\Sigma_0}\Big[{\frac{1}{\gamma}(X_T^{\pi})^{\gamma}}\Big]\leq \frac{x_0^\gamma}{\gamma}\expo\Big(\int_0^T \gamma r_s + \operatorname{Tr}[f(\psi)(T-s)\Sigma_0 +\psi(T-s)NN^{\top}]ds \Big),
\]

which completes the proof. \qed

\section{Martingale property of stochastic exponentials}

In this section we proof the martingale property of the stochastic exponentials appearing in Section \ref{model} and Section \ref{Wishart}. The following lemma is an adaption of \cite[Lemma 7.3]{AJ19} to the multivariate case.

\begin{lem}\label{A1}
Let us denote
\[
M_t:=\expo\Big(\int_0^t \operatorname{Tr}(A_s\Sigma_s^{1/2}dW_s^{\sigma})-\frac{1}{2}\int_0^t\lVert A_s\Sigma_s^{1/2}\rVert^2 ds\Big)
\]
where $(A_t)_{t\in[0,T]}$ is a deterministic process with values in $\R^{d\times d}$ and bounded by $A^*\in\R^{d\times d}$. Then $(M_t)_{t\in[0,T]}$ is a martingale.
\end{lem}    
    
\emph{Proof:}
The process $M_t$ is a stochastic exponential of the form $M_t=\mathcal{E}(\int_0^t \operatorname{Tr}(A_s\Sigma_s^{1/2}dW_s^{\sigma}))$.
Since $M$ is a non-negative local martingale, $M$ is a supermartingale by Fatou's lemma. Thus, in order to show that $M$ is a true martingale, it suffices to show that $\EX[M_T]=1$ for any $T>0$.
For a fixed $T>0$ we define the stopping times 
\[
{\tau_n:=\inf\{{t\geq 0}:\exists 1\leq i,j \leq d: \abs{\Sigma_t^{ij}} >n \}.}
\]
The process $M^{\tau_n}=M_{\tau_n\wedge .}$ is a uniformly integrable martingale for each $n$, since the Novikov condition is fulfilled due to the boundedness of $A$. Thus we get
\[
1=M_0^{\tau_n}=\EX_{\mathbb{P}}[M_T^{\tau_n}]=\EX_{\mathbb{P}}[M_T\vb{1}_{\tau_n\geq T}]+\EX_{\mathbb{P}}[{M_{\tau_n}}\vb{1}_{\tau_n< T}].
\]
By the theorem of dominated convergence, $\EX_{\mathbb{P}}[M_T\vb{1}_{\tau_n\geq T}]\rightarrow \EX_{\mathbb{P}}[M_T]$ and thus, in order to show that $\EX_{\mathcal{P}}[M_T]=1$, it is sufficient to prove that 
\begin{equation*}
    \EX_{\mathbb{P}}[{M_{\tau_n}}\vb{1}_{\tau_n< T}]\rightarrow 0,\text{ as }n\rightarrow\infty.
\end{equation*}
Since $M^{\tau_n}$ is a martingale, we can define probability measures $\mathbb{Q}^n$ with Radon-Nikodym densities
\[
\frac{d\mathbb{Q}^n}{d\mathbb{P}}=M_{\tau_n}.
\]
By Girsanov's theorem the process $W^{\sigma,n}$ defined by
\[
dW_t^{\sigma, n}=dW_t^{\sigma} + \vb{1}_{t\leq \tau_n}\Sigma_t^{1/2}A_t^{\top}ds
\]
is a $d$-dimensional Brownian motions under the measure $\mathbb{Q}^n$. Furthermore, under $\mathbb{Q}^n$ we have 
with
\begin{equation*}
    \begin{aligned}
        \Sigma_t^n {} &=\Sigma_0+\int_0^t K(t-s)(NN^{\top}+M\Sigma_s+\Sigma_s M^{\top}+\vb{1}_{s\leq \tau_n}\Sigma_s A_s^{\top}Q+\vb{1}_{s\leq \tau_n}Q^{\top} A_s\Sigma_s)ds \\
       &\qquad +\int_0^t K(t-s)(\Sigma_s^{1/2}dW_s^{\sigma,n}Q+Q^{\top}(dW_s^{\sigma,n})^{\top}\Sigma_s^{1/2}).
    \end{aligned}
\end{equation*}
Using the vectorization operator from Section \ref{Existence}, the above equation can be written as
\begin{equation*}
    \begin{aligned}
        \vecop(\Sigma^n) {} &=\vecop(\Sigma_0)+ (I\otimes K)\ast[\vecop(NN^{\top})\\
        &\qquad\qquad\qquad +(I\otimes (M+\vb{1}_{\leq \tau_n}Q^{\top} A)+(M+\vb{1}_{\leq \tau_n}Q^{\top} A)\otimes I)\vecop(\Sigma)]\\
       &\qquad +(I\otimes K)\ast  [Q^{\top}\otimes\Sigma^{1/2}\vecop(dW^{\sigma,n})+\Sigma^{1/2}\otimes Q\vecop((dW^{\sigma,n})^{\top})]\\
       &= (I\otimes K)\ast[b(\vecop(\Sigma))dt + c(\vecop(\Sigma))d\vecop(W^{\sigma,n})].
    \end{aligned}
\end{equation*}
where 
\[
b(\vecop(\Sigma)):=\vecop(NN^{\top}) +(I\otimes (M+\vb{1}_{\leq \tau_n}Q^{\top} A)+(M+\vb{1}_{\leq \tau_n}Q^{\top} A)\otimes I)\vecop(\Sigma)
\]
and 
\[
c(\vecop(\Sigma)):=Q^{\top}\otimes \Sigma^{1/2}+f(\Sigma^{1/2}\otimes Q^{\top})
\]
{and $f: \R^{d\times d} \rightarrow \R^{d\times d}$ swaps entries such that $f(\Sigma^{1/2}\otimes Q^{\top})d\operatorname{vec}(W^{\sigma})=\Sigma^{1/2}\otimes Q^{\top}d\operatorname{vec}((W^{\sigma})^{\top})$}.
Using similar arguments as in the proof of Theorem \ref{ExVol}, one can show that the drift and the diffusion term of the above equation fulfill the linear growth condition \cite[condition (3.1)]{AJ19}, i.e. we have 
\[
\abs{b(\vecop(\Sigma))}\vee \abs{c(\vecop(\Sigma))}\leq c_{LG}(1+\abs{\vecop(\Sigma})).
\]
Note that the argument for the drift only works if the matrix $A_t$ is bounded. Choose $p>2$ sufficiently large that $\kappa_i/2-1/p>0$ for $\kappa_i$ defined in \eqref{kernels}. An application of \cite[Lemma 3.1]{AJ19} yields the moment bound
\[
\sup_{t\leq T}\EX[\abs{\vecop(\Sigma_t)}^p]\leq c
\]
for some constant $c$ independent of $n$. The $0$-Hölder seminorm of a function $f$ is defined as
\[
\abs{f}_{C^{0,0}(0,T)}=\sup_{0\leq s<t\leq T}\abs{f(t)-f(s)}.
\]

$\textbf{Claim:}$ For some constant $C$ independent of $n$ the following inequality holds:
\[
\sum\limits_{1\leq i,j\leq d}\abs{\Sigma^{ij}}_{C^{0,0}}^p\leq C\abs{\vecop(\Sigma)}_{C^{0,0}}^p.
\]
\emph{Proof of Claim:} We have to show that
\[
\sum\limits_{1\leq i,j\leq d}(\sup_{0\leq s<t\leq T}\abs{\Sigma_t^{ij}-\Sigma_s^{ij}})^p\leq C (\sup_{0\leq s<t\leq T}(\sum\limits_{1\leq i,j\leq d}\abs{\Sigma_t^{ij}-\Sigma_s^{ij}}^2)^{1/2})^p
\]
or equivalently
\[
(\sum\limits_{1\leq i,j\leq d}(\sup_{0\leq s<t\leq T}\abs{\Sigma_t^{ij}-\Sigma_s^{ij}})^p)^{1/p}\leq C \sup_{0\leq s<t\leq T}(\sum\limits_{1\leq i,j\leq d}\abs{\Sigma_t^{ij}-\Sigma_s^{ij}}^2)^{1/2}.
\]
Since for a vector $x$ we have $\lVert x\rVert_{l^q}\leq \lVert x\rVert_{l^p}$ for $1\leq p<q\leq\infty$, it holds that the left-hand side is bounded by
\[
(\sum\limits_{1\leq i,j\leq d}(\sup_{0\leq s<t\leq T}\abs{\Sigma_t^{ij}-\Sigma_s^{ij}})^p)^{1/p}\leq (\sum\limits_{1\leq i,j\leq d}(\sup_{0\leq s<t\leq T}\abs{\Sigma_t^{ij}-\Sigma_s^{ij}})^2)^{1/2}.
\]
Using the fact that the square of the supremum of a set of non-negative numbers equals the supremum of the squares, it remains to show that
\[
\sum\limits_{1\leq i,j\leq d}\sup_{0\leq s<t\leq T}\abs{\Sigma_t^{ij}-\Sigma_s^{ij}}^2\leq C \sup_{0\leq s<t\leq T}\sum\limits_{1\leq i,j\leq d}\abs{\Sigma_t^{ij}-\Sigma_s^{ij}}^2.
\]
Clearly
\[
\sum\limits_{1\leq i,j\leq d}\sup_{0\leq s<t\leq T}\abs{\Sigma_t^{ij}-\Sigma_s^{ij}}^2\leq d^2 \max_{1\leq i,j\leq d}(\sup_{0\leq s<t\leq T}\abs{\Sigma_t^{ij}-\Sigma_s^{ij}}^2).
\]
and since
\[
\forall 1\leq i,j\leq d,\, \forall 0\leq s<t\leq T: \abs{\Sigma_t^{ij}-\Sigma_s^{ij}}^2\leq \sum\limits_{1\leq i,j\leq d} \abs{\Sigma_t^{ij}-\Sigma_s^{ij}}^2
\]
we get 
\[
\max_{1\leq i,j\leq d}(\sup_{0\leq s<t\leq T}\abs{\Sigma_t^{ij}-\Sigma_s^{ij}}^2)\leq \sup_{0\leq s<t\leq T}\sum\limits_{1\leq i,j\leq d}\abs{\Sigma_t^{ij}-\Sigma_s^{ij}}^2.
\]
This completes the proof of the claim. \qed

We now show that $\EX_{\mathbb{P}}[M_T\vb{1}_{\tau_n< T}]\rightarrow 0$ as $n\rightarrow\infty$:
\begin{equation*}
    \begin{aligned}
       \EX_{\mathbb{P}}[{M_{\tau_n}}\vb{1}_{\tau_n< T}] {}
       &= \Q^n(\tau_n<T)=\Q^n(\exists 1\leq i,j\leq d: \sup_{t\leq T}\abs{\Sigma_t^{ij}}>n)\\
       &\leq \sum\limits_{1\leq i,j \leq d} \Q^n(\sup_{t\leq T}\abs{\Sigma_t^{ij}}>n) \leq \sum\limits_{1\leq i,j \leq d} \Q^n(\abs{\Sigma_0^{ij}}+\abs{\Sigma^{ij}}_{C^{0,0}(0,T)}>n)\\
       &\leq \sum\limits_{1\leq i,j\leq d} (\frac{1}{n-\abs{\Sigma_0^{ij}}})^p\EX_{\Q^n}[\abs{\Sigma^{ij}}_{C^{0,0}(0,T)}^p]\text{ (Markov inequality)}\\
       &\leq (\frac{1}{n-\abs{\Sigma_0^{\operatorname{max}}}})^p \EX_{\Q^n}[\sum\limits_{1\leq i,j\leq d}\abs{\Sigma^{ij}}_{C^{0,0}(0,T)}^p]\\
       &\leq C(\frac{1}{n-\abs{\Sigma_0^{\operatorname{max}}}})^p \EX_{\Q^n}[\abs{\vecop(\Sigma)}_{C^{0,0}(0,T)}^p]\text{ (Claim)}\\
       &\leq C (\frac{1}{n-\abs{\Sigma_0^{\operatorname{max}}}})^p\sup_{t\leq T}\EX_{\Q^n}[\abs{b(\vecop(\Sigma_t))}^p+\abs{c(\vecop(\Sigma_t))}^p]\text{ (\cite[Lemma 2.4]{AJ19})}\\
       &\leq C (\frac{1}{n-\abs{\Sigma_0^{\operatorname{max}}}})^p\sup_{t\leq T}\EX_{\Q^n}[(1+\abs{\vecop(\Sigma_t)})^p]\text{ (Growth condition)}\\
       &\leq C (\frac{1}{n-\abs{\Sigma_0^{\operatorname{max}}}})^p\sup_{t\leq T}\EX_{\Q^n}[1+\abs{\vecop(\Sigma_t)}^p]\\
       &\leq C (\frac{1}{n-\abs{\Sigma_0^{\operatorname{max}}}})^p\text{ (\cite[Lemma 3.1]{AJ19})}.
    \end{aligned}
\end{equation*}
Since the constant $C$ is independent of $n$ we finally obtain
\[
\EX_{\mathbb{P}}[{M_{\tau_n}}\vb{1}_{\tau_n< T}]\leq C (\frac{1}{n-\abs{\Sigma_0^{\operatorname{max}}}})^p\xrightarrow{n\rightarrow\infty} 0,
\]
which completes the proof. \qed\\

\begin{lem}\label{A2}
Let us denote
\[
M_t:=\expo\Big(\int_t^T \operatorname{Tr}(A_s^{\top}\Sigma_s^{1/2}dW_s^{\sigma})-\frac{1}{2}\int_t^T\lVert A_s^{\top}\Sigma_s^{1/2}\rVert^2 ds\Big)
\]
where $(A_t)_{t\in[0,T]}$ is a deterministic process with values in $\R^{d}$ and bounded by $A^*\in\R^{d}$. Then $(M_t)_{t\in[0,T]}$ is a martingale.
\end{lem}

\emph{Proof:} This is a direct consequence of Lemma~\ref{A1} and \cite[Proposition A.2.]{BL13}.\qed\\

The next lemma is an enhancement of \cite[Appendix C]{PJM21}, the proof follows similar arguments.  

\begin{lem}\label{extended AJ lemma}
Let $W_1$, $W_2$ be two independent $d$-dimensional brownian motions and for $1\leq i\leq d$ let $g_{1i},  g_{2i}\in L^{\infty}(\R^+,\R)$. Then the local martingale 
\[
Z_t=\mathcal{E}(\int_0^t\sum\limits_{i=1}^d g_{1i}(s)\sqrt{V_s^i}dW_{1s}^i+\int_0^t\sum\limits_{i=1}^d g_{2i}(s)\sqrt{V_s^i}dW_{2s}^i)
\]
is a true martingale.
\end{lem}    
    
\emph{Proof:}
Set $U=\int_0^. V_s ds$. Then by the stochastic Fubini theorem we get
\[
U_t^i=\int_0^t v_0^i(s)ds +\int_0^t K_i(t-s)Z_s^i ds
\]
with
\[
Z_t^i=\int_0^t (DV_s)_i ds+\int_0^t \nu_i\sqrt{V_s^i}dB_s^i.
\]
Since $M$ is a non-negative local martingale, $M$ is a supermartingale by Fatou's lemma. Thus, in order to show that $M$ is a true martingale, it suffices to show that $\EX[M_T]=1$ for any $T>0$.
For a fixed $T>0$ we define the stopping times 
\[
\tau_n:=\inf\{{t\geq 0}:\exists 1\leq i\leq d: \int_0^t V_s^i ds>n \}.
\]
The stopped process $M^{\tau_n}=M_{\tau_n\wedge .}$ is a uniformly integrable martingale for each $n$, since the Novikov condition is fulfilled due to the boundedness of $g_1$ and $g_2$. Thus we get
\[
1=M_0^{\tau_n}=\EX_{\mathbb{P}}[M_T^{\tau_n}]=\EX_{\mathbb{P}}[M_T\vb{1}_{\tau_n\geq T}]+\EX_{\mathbb{P}}[M_T\vb{1}_{\tau_n< T}].
\]
By the theorem of dominated convergence, $\EX_{\mathbb{P}}[M_T\vb{1}_{\tau_n\geq T}]\rightarrow \EX_{\mathbb{P}}[M_T]$ and thus, in order to show that $\EX_{\mathcal{P}}[M_T]=1$, it is sufficient to prove that 
\begin{equation}
    \EX_{\mathbb{P}}[M_T\vb{1}_{\tau_n< T}]\rightarrow 0,\text{ as }n\rightarrow\infty.
\end{equation}
Since $M^{\tau_n}$ is a martingale, we can define probability measures $\mathbb{Q}^n$ with Radon-Nikodym densities
\[
\frac{d\mathbb{Q}^n}{d\mathbb{P}}=M_{\tau_n}.
\]
By Girsanov's theorem the processes $W_1^n$ and $W_2^n$ defined by
\[
W_1^{n,i}=W_1^i +\int_0^. \vb{1}_{s\leq \tau_n} g_{1,i}(s)\sqrt{V_s^i}ds,\, 1\leq i\leq d,
\]
\[
W_2^{n,i}=W_2^i +\int_0^. \vb{1}_{s\leq \tau_n} g_{2,i}(s)\sqrt{V_s^i}ds,\, 1\leq i\leq d,
\]
are $d$-dimensional Brownian motions under the measure $\mathbb{Q}^n$. Furthermore, under $\mathbb{Q}^n$ we have 
\[
U_t^i=\int_0^t v_0^i(s)ds+\int_0^t K_i(t-s)Z_s^{n,i}ds
\]
\begin{equation*}
    \begin{aligned}
       Z_t^{n,i} {} &=\int_0^t ((DV_s)_i-\vb{1}_{s\leq \tau_n}\rho_i\nu_ig_{1,i}(s)V_s^i-\vb{1}_{s\leq \tau_n})\sqrt{1-\rho_i^2}\nu_i g_{2,i}(s)V_s^i)ds \\
       &\qquad +\int_0^t \nu_i\sqrt{V_s^i}(\rho_i dW_{1s}^{n,i}+\sqrt{1-\rho_i^2}dW_{2s}^{n,i})
    \end{aligned}
\end{equation*}
and under $\mathbb{Q}^n$, the drift of $Z^n$ satisfies a linear growth condition in $U$ for some constant $\kappa_L$ independent of $n$. Therefore an application of the generalized Grönwall inequality (cf. \cite[Lemma 3.1]{AJ21}) yields the moment bound
\[
\EX_{\mathbb{Q}^n}[\abs{U_T}^2]\leq \eta(\kappa_L,T,K,v_0),
\]
where $\eta(\kappa_L,T,K,v_0)$ does not depend on $n$. An application of Chebyshev's inequality yields
\begin{equation*}
    \begin{aligned}
       \EX_{\mathbb{P}}[M_T\vb{1}_{\tau_n< T}] {}
       &= \Q^n(\tau_n<T)\\
       &\leq \sum\limits_{i=1}^d \Q^n(U_T^i>n)\\
       &\leq \sum\limits_{i=1}^d \frac{1}{n^2}\EX_{\Q^n}[\abs{U_T^i}^2]\\
       &= \frac{1}{n^2}\EX_{\Q^n}[\abs{U_T}^2]\\
       &\leq \frac{1}{n^2} \eta(\kappa_L,T,K,v_0)\xrightarrow{n\rightarrow\infty}0.
    \end{aligned}
\end{equation*}

This completes the proof. \qed

\section{Dynamics of the process $M$}
We derive the dynamics of the process $M$ appearing in the martingale distortion approach using Itô's formula (cf. \cite[Theorem 3.2]{HW21}).  

\begin{lem}\label{dynamics M}
The process $M_t$ defined in \eqref{M} has dynamics
\begin{equation*}
    \begin{aligned}
dM_t {}
    &= M_t[-\gamma r_t-\frac{\gamma}{2(1-\gamma)}\theta\Theta V_t]dt\\
    &\qquad + M_t c\psi(T-t) N \sqrt{\diag{(V_t)}}P_1 d{W}_{1t} + M_t c\psi(T-t) N \sqrt{\diag{(V_t)}}P_2 d{W}_{2t}\\
    &\qquad - \frac{\gamma}{1-\gamma}M_t c\psi(T-t)N P_1\diag(V_t)\theta^{\top}-\frac{\gamma}{2(1-\gamma)} M_t \lVert c\psi(T-t) N \sqrt{\diag{(V_t)}}P_1 \rVert_2^2.
    \end{aligned}
\end{equation*}
\end{lem}

\emph{Proof:}
Let $Z_t =\int_t^T[\gamma r_s +\frac{\gamma}{2(1-\gamma)}\theta\Theta\xi_t(s)+\frac{c}{2}A(\psi)(T-s)\xi_t(s)]ds$. Then $M_t=e^{Z_t}$. Applying Itô's lemma to $\xi_t(s)$ yields
\[
d\xi_t(s)=R_{\Lambda}(s-t)\Lambda^{-1}N\sqrt{\diag(V_t)}d\tilde{B}_t
\]
by \cite[Lemma 4.2]{AJ19}. Define $A(\psi):=\psi N^2\Psi$.
Then

\begin{equation*}
    \begin{aligned}
dZ_t {}
    &= [-\gamma r_t-\frac{\gamma}{2(1-\gamma)}\theta\Theta V_t-\frac{c}{2}A(\psi)(T-t)V_t]dt\\
    &\quad+ \int_t^T\frac{\gamma}{2(1-\gamma)}\theta\Theta R_{\Lambda}(s-t)\Lambda^{-1}N\sqrt{\diag{(V_t)}}d\tilde{B}_t ds\\
    &\quad+ \int_t^T \frac{c}{2}A(\psi)(T-s)R_{\Lambda}(s-t)\Lambda^{-1}N\sqrt{\diag{(V_t)}}d\tilde{B}_t ds\\
    &= [-\gamma r_t-\frac{\gamma}{2(1-\gamma)}\theta\Theta V_t-\frac{c}{2}A(\psi)(T-t)V_t]dt\\
    &\quad+ \int_t^T\frac{\gamma}{2(1-\gamma)}\theta\Theta R_{\Lambda}(s-t)\Lambda^{-1}N ds\sqrt{\diag{(V_t)}}d\tilde{B}_t \\
    &\quad+ \int_t^T \frac{c}{2}A(\psi)(T-s)R_{\Lambda}(s-t)\Lambda^{-1}N ds\sqrt{\diag{(V_t)}}d\tilde{B}_t\\
    &= [-\gamma r_t-\frac{\gamma}{2(1-\gamma)}\theta\Theta V_t-\frac{c}{2}A(\psi)(T-t)V_t]dt\\
    &\quad+ \int_t^T[\frac{c}{2}A(\psi)(T-s)+\frac{\gamma}{2(1-\gamma)}\theta\Theta]R_{\Lambda}(s-t)\Lambda^{-1} ds N \sqrt{\diag{(V_t)}}d\tilde{B}_t.
    \end{aligned}
\end{equation*}
Here, for the second equality, we used the stochastic Fubini theorem from \cite{Ver12}.
Next, we show that
\[
\int_t^T[\frac{c}{2}A(\psi)(T-s)+\frac{\gamma}{2(1-\gamma)}\theta\Theta]R_{\Lambda}(s-t)\Lambda^{-1} ds =c\psi(T-t).
\]
We have
\begin{equation*}
    \begin{aligned}
    {}
    & \int_t^T[\frac{c}{2}A(\psi)(T-s)+\frac{\gamma}{2(1-\gamma)}\theta\Theta]R_{\Lambda}(s-t)\Lambda^{-1} ds -c\psi(T-t)\\
    &= [\frac{c}{2}A(\psi)+\frac{\gamma}{2(1-\gamma)}\theta\Theta]\ast (R_{\Lambda} \Lambda^{-1})(T-t)\\
    &\qquad\qquad\qquad\qquad\qquad\qquad\quad -[\frac{c}{2}A(\psi)-\psi\Lambda+\frac{\gamma}{2(1-\gamma)}\theta\Theta]\ast(cK)(T-t)\\
    &= [\frac{c}{2}A(\psi)+\frac{\gamma}{2(1-\gamma)}\theta\Theta]\ast (R_{\Lambda}\Lambda^{-1}-K)(T-t)+c(\psi \Lambda)\ast K(T-t)\\
    &= [\frac{c}{2}A(\psi)-c\psi\Lambda+c\psi\Lambda+\frac{\gamma}{2(1-\gamma)}\theta\Theta]\ast (R_{\Lambda}\Lambda^{-1}-K)(T-t)+c(\psi \Lambda)\ast K(T-t)\\
    &= [\frac{c}{2}A(\psi)-c\psi\Lambda+\frac{\gamma}{2(1-\gamma)}\theta\Theta]\ast((-K\Lambda\ast R_{\Lambda})\Lambda^{-1})\\
    &\qquad\qquad\qquad\qquad\qquad\qquad\quad + c(\psi\Lambda)\ast(-R_{\Lambda}\ast K)(T-t)+c(\psi \Lambda)\ast K(T-t)\\
    &= (-c(\psi\Lambda)\ast R_{\Lambda}\Lambda^{-1})(T-t)-c(\psi\Lambda)\ast(R_{\Lambda}\ast K)(T-t)+c(\psi \Lambda)\ast K(T-t)\\
    &= c(\psi\Lambda)\ast[K-R_{\Lambda}\ast K-R_{\Lambda}\Lambda^{-1}](T-t) = 0.
    \end{aligned}
\end{equation*}
Here the last equality holds, since
\begin{equation*}
    \begin{aligned}
    K-R_{\Lambda}\ast K {}
    &= (K-R_{\Lambda}\ast K)\Lambda\Lambda^{-1}=(K\Lambda-(R_{\Lambda}\ast K)\Lambda)\Lambda^{-1}\\
    &= (K\Lambda-R_{\Lambda}\ast(K\Lambda))\Lambda^{-1}=(K\Lambda-(K\Lambda-R_{\Lambda}))\Lambda^{-1}=R_{\Lambda}\Lambda^{-1}.
    \end{aligned}
\end{equation*}    
Thus we get 

\begin{equation*}
    \begin{aligned}
dZ_t {}
    &= [-\gamma r_t-\frac{\gamma}{2(1-\gamma)}\theta\Theta V_t-\frac{c}{2}A(\psi)(T-t)V_t]dt + c\psi(T-t) N \sqrt{\diag{(V_t)}}d\tilde{B}_t\\
    &= [-\gamma r_t-\frac{\gamma}{2(1-\gamma)}\theta\Theta V_t-\frac{c}{2}A(\psi)(T-t)V_t]dt\\
    &\qquad\qquad + c\psi(T-t) N \sqrt{\diag{(V_t)}}P_1 d\tilde{W}_{1t} + c\psi(T-t) N \sqrt{\diag{(V_t)}}P_2 d{W}_{2t}
    \end{aligned}
\end{equation*}

where $P_1=\diag{(\rho)}$ and $P_2=\diag{(\sqrt{1-\rho^2})}$.

Since $M_t=e^{Z_t}$, by Itô's formula we obtain $dM_t=M_tdZ_t+\frac{1}{2}M_td\langle Z_t,Z_t\rangle$, i.e

\begin{equation*}
    \begin{aligned}
dM_t {}
    &= M_t[-\gamma r_t-\frac{\gamma}{2(1-\gamma)}\theta\Theta V_t-\frac{c}{2}A(\psi)(T-t)V_t]dt\\
    &\qquad + M_t c\psi(T-t) N \sqrt{\diag{(V_t)}}P_1 d\tilde{W}_{1t} + M_t c\psi(T-t) N \sqrt{\diag{(V_t)}}P_2 d{W}_{2t}\\
    &\qquad + \frac{1}{2}M_t \lVert c\psi(T-t) N \sqrt{\diag{(V_t)}}P_1 \rVert_2^2 + \frac{1}{2}M_t \lVert c\psi(T-t) N \sqrt{\diag{(V_t)}}P_2 \rVert_2^2\\
    &= M_t[-\gamma r_t-\frac{\gamma}{2(1-\gamma)}\theta\Theta V_t]dt\\
    &\qquad + M_t c\psi(T-t) N \sqrt{\diag{(V_t)}}P_1 d{W}_{1t} + M_t c\psi(T-t) N \sqrt{\diag{(V_t)}}P_2 d{W}_{2t}\\
    &\qquad - \frac{\gamma}{1-\gamma}M_t c\psi(T-t)N P_1\diag(V_t)\theta^{\top}-\frac{\gamma}{2(1-\gamma)} M_t \lVert c\psi(T-t) N \sqrt{\diag{(V_t)}}P_1 \rVert_2^2.
    \end{aligned}
\end{equation*}

The last equality holds, since in the degenerate correlation case we have the following identity:

\begin{equation*}
    \begin{aligned}
    {}
    & -M_t\frac{c}{2}A(\psi)(T-t)V_tdt+ \frac{1}{2}M_t \lVert c\psi(T-t) N \sqrt{\diag{(V_t)}}P_2 \rVert_2^2\\
    &= -\frac{c}{2}M_t\sum\limits_{i=1}^d \psi_i^2(T-t)\nu_i^2 V_i dt + \frac{1}{2}M_t \sum\limits_{i=1}^d c^2\psi_i^2(T-t)\nu_i^2 V_i (1-\rho^2)dt\\
    &= (-\frac{1}{c}+(1-\rho^2))\frac{1}{2}M_t \sum\limits_{i=1}^d c^2\psi_i^2(T-t)\nu_i^2 V_i dt\\
    &= -\frac{\rho^2}{1-\gamma}\frac{1}{2}M_t\sum\limits_{i=1}^d c^2\psi_i^2(T-t)\nu_i^2 V_i dt
    = -\frac{1}{1-\gamma}\frac{1}{2}M_t\sum\limits_{i=1}^d c^2\psi_i^2(T-t)\nu_i^2 V_i \rho^2 dt\\
    & = -\frac{1}{1-\gamma}\frac{1}{2} M_t \lVert c\psi(T-t) N \sqrt{\diag{(V_t)}}P_1 \rVert_2^2. 
    \end{aligned}
\end{equation*}

\qed

\section*{Acknowledgements} F. Aichinger and S. Desmettre are supported by the Austrian Science Fund (FWF) project  F5507-N26, which is part of the Special Research Program \textit{Quasi-Monte Carlo Methods: Theory and Applications}. We also thank Bingyan Han for useful discussions concerning Lemma~\ref{extended AJ lemma}.

\end{document}